\documentclass[12pt,twoside]{article}
\usepackage[hmargin=0.8in,vmargin=0.9in]{geometry}
\usepackage{fancyhdr}
\usepackage{graphicx}
\usepackage{subfigure}
\usepackage{amssymb}
\usepackage{amsmath}
\usepackage{amsfonts}
\usepackage{theorem}
\usepackage{mathrsfs}
\usepackage{mathtools}
\usepackage{empheq}
\usepackage{bm}
\usepackage{bbm}
\usepackage{color}
\usepackage{setspace}
\usepackage{epstopdf}
\usepackage{float}
\usepackage{picinpar}
\usepackage{extarrows}
\usepackage{multirow}
\usepackage{hyperref}

\usepackage{cite}
\usepackage{nicefrac}
\usepackage{booktabs} 
\usepackage{array} 
\usepackage{paralist} 
\usepackage{verbatim} 
\usepackage{subfigure} 
\usepackage{authblk}
\usepackage{cleveref}

\allowdisplaybreaks[1]
\usepackage[normalem]{ulem}

\newtheorem{thm}{Theorem}[section]

\newtheorem{rmk}[thm]{Remark}

{\theorembodyfont{\rmfamily}}

\newcommand\eref[1]{(\ref{#1})}

\newcommand*\xbar[1]{%
  \hbox{%
    \vbox{%
      \hrule height 0.5pt 
      \kern0.4ex
      \hbox{%
        \kern-0.05em
        \ensuremath{#1}%
        \kern-0.00em
      }%
    }%
  }%
}

\allowdisplaybreaks[1]

\def\hf {\frac{1}{2}}

\newcommand{\jph}{{j+\frac{1}{2}}}

\newcommand{\kph}{{k+\frac{1}{2}}}
\newcommand{\kmh}{{k-\frac{1}{2}}}

\newcommand{\dx}{\Delta x}

\newcommand{\dt}{\Delta t}

\numberwithin{equation}{section}
\numberwithin{figure}{section}
\numberwithin{table}{section}

\graphicspath{{figures/}{./}}

\title{Flux Globalization Based Well-Balanced Path-Conservative Central-Upwind Schemes for Shallow Water Linearized Moment Equations}
\author{Yangyang Cao\thanks{MSU-BIT-SMBU Joint Research Center of Applied Mathematics, Shenzhen MSU-BIT University, Shenzhen, 518172, China;
{\tt caoyangyang@smbu.edu.cn}},~Qian Huang\thanks{Institute of Applied Analysis and Numerical Simulation, University of Stuttgart,
70569 Stuttgart, Germany; {\tt qian.huang@mathematik.uni-stuttgart.de}},~Julian Koellermeier\thanks{Department of Mathematics, Computer
Science and Statistics, Ghent University, 9000 Gent, Belgium; {\tt julian.koellermeier@ugent.be} and Bernoulli Institute, University of
Groningen, 9747 AG Groningen, The Netherlands; {\tt j.koellermeier@rug.nl}},~Alexander Kurganov\thanks{Department of Mathematics and
Shenzhen International Center for Mathematics, Southern University of Science and Technology, Shenzhen, 518055, China; {\tt
alexander@sustech.edu.cn}},~and Yongle Liu\thanks{Institute of Mathematics, University of Z\"{u}rich, Z\"{u}rich, 8057, Switzerland;
{\tt liuyl2017@mail.sustech.edu.cn}}}

\begin{document}
\date{}
\maketitle
\begin{abstract}
We develop second-order path-conservative central-upwind (PCCU) schemes for the  hyperbolic shallow water linearized moment equations
(HSWLME), which are an extension of standard depth-averaged models for free-surface flows. The proposed PCCU schemes are constructed via
flux globalization strategies adapted to the nonconservative form via a path-conservative finite-volume method. The resulting scheme is
well-balanced (WB) in the sense that it is capable of exactly preserving physically relevant steady states including moving-water ones. We
validate the proposed scheme  on several benchmarks, including smooth solutions, small perturbation of steady states, and dam-break
scenarios. These results demonstrate that our flux globalization based WB PCCU schemes provide a reliable framework for computing solutions
of shallow water moment models with nonlinear and nonconservative features.
\end{abstract}



\section{Introduction}\label{sec1}
The accurate computation of free-surface flows depends on the combination of mathematical models and tailored numerical methods. The most
simple one-dimensional shallow water equations (SWE) result from depth-averaging the incompressible Navier-Stokes equations, but only take
into account the water depth $h(t,x)$ and the depth-averaged velocity $u(t,x)$, therefore not modeling the dynamical evolution of variations
from a constant velocity profile. Various models improve upon the standard SWE, for example, by considering multiple layers
\cite{AB2007,ABPelS,ABPerS,BFGN,FKC} or non-hydrostatic effects \cite{BBBD,BMSS,FPPS,YKC}. Recently, families of so-called moment models
that include evolution equations for a polynomial expansion of a velocity profile have been derived. They all assume an expansion of the
horizontal velocity $u(t,x,z)$ in vertical direction $z$, that is, $u(t,x,z)=u(t,x)+\sum_{i=1}^N\alpha_i\phi_i(z)$, with depth-averaged
velocity $u$, basis coefficients $\alpha_i$, scaled Legendre polynomials $\phi_i$, and expansion order $N$; see \cite{KT_SWM19}.

This paper focuses on the development of highly accurate and robust numerical methods for the arbitrary order hyperbolic shallow water
linearized moment equations (HSWLME) that take into account nonflat bottom topography and friction terms. In the general $N$-moment case,
the HSWLME read as (see, e.g., \cite{KP_SWM})
\begin{equation}
\bm U_t+\bm F(\bm U)_x=Q(\bm U)\bm U_x+\bm S(\bm U)+\bm P(\bm U),
\label{1.1}
\end{equation}
where
\begin{equation}
\begin{aligned}
&\bm U=(h,hu,h\alpha_1,\ldots,h\alpha_N)^\top,\quad\bm S(\bm U)=(0,-gh\partial_x Z,0,\ldots,0)^\top,\quad Q={\rm diag}(0,0,u,\ldots,u),\\
&\bm F(\bm U)=\Big(hu,hu_2+\hf gh^2+\frac{1}{3}h\alpha_1^2+\ldots+\frac{1}{2N+1}h\alpha_N^2,2hu\alpha_1,\ldots,2hu\alpha_N\Big)^\top,
\label{1.2}
\end{aligned}
\end{equation}
and $\bm P(\bm U)$ is a friction term. In \eref{1.2}, $x$ is a spatial variable, $t$ is time, $h=h(x,t)$ is the water depth, $u=u(x,t)$ is
the mean horizontal water velocity, $g$ is an acceleration due to gravity, $\alpha_i=\alpha_i(x,t)$, $i=1,\ldots,N$ are coefficients of the
velocity profile, and $Z(x)$ is the bottom topography. Finally, the friction term  $\bm P=(0,P_0,P_1,\ldots,P_N)^\top$ is given by
\begin{equation}
\begin{aligned}
&P_0=-\frac{\nu}{\lambda}\Big(u+\sum\limits_{j=1}^N\alpha_j\Big),\\
&P_i=-\frac{\nu}{\lambda}(2i+1)\Big(u+\sum\limits_{j=1}^N\alpha_j\Big)-4\frac{\nu}{h}(2i+1)\sum\limits_{j=1}^Na_{i,j}\alpha_j,\quad
i=1,\ldots,N,
\end{aligned}
\label{1.3}
\end{equation}
where $\nu$ and $\lambda$ are positive coefficients representing the kinematic viscosity and slip length, respectively, and
\begin{equation}
a_{i,j}=\left\{\begin{aligned}
&0&&\mbox{if $i+j$ is even},\\
&\hf\min(i-1,j)(\min(i-1,j)+1)&&\mbox{if $i+j$ is odd}.
\end{aligned}\right.
\label{1.4}
\end{equation}

Shallow water moment models, which take into account horizontal velocity changes over the vertical direction, were first developed in
\cite{KT_SWM19} based on the Legendre polynomial-based expansions of the velocity modeling the deviation from a constant velocity profile. A
further analysis showed that those models suffer from loss of hyperbolicity, which was solved by a hyperbolic regularization strategy in
\cite{KR_2020}, resulting in a globally hyperbolic shallow water moment model. Followed by successful applications and extensions of the
models (see e.g., \cite{Scholz2023,Koellermeier2023,Garres-Diaz2021a,Huang2022,Garres-Diaz2023}), steady states of shallow water moment
models were first investigated in \cite{KP_SWM}, which led to the development of the HSWLME, that can be considered as the simplest shallow
water moment model, in which the nonlinear contributions of the basis coefficients is neglected in the higher-order moment equations, while
keeping the mass and momentum equations unchanged. We refer the reader to \cite{KP_SWM} for the hyperbolicity, steady states, and
Rankine-Hugoniot conditions of the frictionless version of the HSWLME.

Designing a good numerical scheme for \eref{1.1}--\eref{1.2} is a challenging task for a number of reasons. First, the presence of the
nonconservative product term $Q(\bm U)\bm U_x$ makes it impossible to understand weak solutions of \eref{1.1} in the sense of distributions
(instead weak solutions can be understood as Borel measures; see \cite{DLM,LPG2002,LPG2004}). Second, the well-balanced (WB) property is
crucial (and not easy to achieve, especially due to the presence of the friction terms), as it allows to exactly preserve some of the
physically relevant steady-state solutions satisfying
\begin{equation}
\bm F(\bm U)_x-Q(\bm U)\bm U_x-\bm S(\bm U)+\bm P(\bm U)=\bm0,
\label{1.5}
\end{equation}
and hence to capture their small perturbation on practically affordable coarse meshes. Note that the presence of the friction term
$\bm P(\bm U)$ in \eref{1.5} makes the structure of the steady states substantially more complicated than in the frictionless case studied
in \cite{KP_SWM}.

The concept of measure-valued solutions was utilized in \cite{CLMP,Par2006} to develop path-conservative numerical methods, which were then
successfully applied to a variety of nonconservative hyperbolic systems; see, e.g.,
\cite{DHZ,CKM,BDGI,CMP2017,CP2020,Cha,PC2009,PGCCMP,SGBNP} and references therein. Most of these path-conservative schemes are WB, but they
can often only preserve very simple steady states like the ``lake-at-rest'' ones since the WB property is enforced in a rather artificial
way by modifying the numerical viscosity terms.

A more systematic approach to develop WB schemes for hyperbolic systems of balance laws is based on a flux globalization, which was
originally introduced in \cite{CDH2009,DDMA11,GC2001,MGD11}. To apply a flux globalization approach to \eref{1.1}, one has to first rewrite
it in a quasi-conservative form
\begin{equation}
\begin{aligned}
&\bm U_t+\bm K_x=\bm0,
\end{aligned}
\label{1.6}
\end{equation}
where $\bm K(\bm U)$ is the following global flux:
\begin{equation}
\bm K(\bm U)=\bm F(\bm U)-\bm R(\bm U),\quad\bm R(\bm U)=\int\limits_{\widehat x}^x\Big[Q(\bm U(\xi,t))\bm U_\xi(\xi,t)+
\bm S(\bm U(\xi,t))+\bm P(\bm U(\xi,t))\Big]{\rm d}\xi,
\label{1.7}
\end{equation}
where $\widehat x$ is an arbitrary number. Steady states can then be written as $\bm K(\bm U)\equiv{\bf Const}$ and according to the flux
globalization method introduced in \cite{CHO_18}, one can design a WB scheme for \eref{1.6}--\eref{1.7} by reconstructing  $\bm K$ instead
of $\bm U$ and by applying a Riemann-problem-solver-free central-upwind (CU) scheme to \eref{1.6}. This method was successfully applied to
several hyperbolic systems of balanced laws, in, e.g., \cite{CCHKW_18,KLZ,CKLHy,CCKOT_18,ChuKur}.

These flux globalization based schemes take advantage of the CU numerical fluxes, which were originally introduced in \cite{KLin,KNP,KTcl}
for general multidimensional nonlinear hyperbolic (systems of) PDEs and can be directly implemented for approximating global fluxes. The
applicability of the aforementioned flux globalization based CU schemes, however, is limited to the case when the global terms
$\bm R(\bm U)$ in \eref{1.7} are continuous, which is not the case when the nonconservative term $Q(\bm U)\bm U_x$ is present. In this case,
the integral in \eref{1.7} should be computed using a path-conservative integration technique as it was done in \cite{KLX_21}, where new
flux globalization based WB path-conservative CU (PCCU) schemes were introduced. In these schemes, we first identify the equilibrium
variables $\bm E(\bm U)$, then reconstruct them, use the reconstructed values of $\bm E$  to obtain the corresponding values of $\bm U$, and
evaluate the global CU numerical fluxes. The new flux globalization based WB PCCU schemes are capable of preserving a much wider range of
steady states, as was demonstrated in \cite{CKL23,CKLX_22,CKLZ,CKRZ_MHD,CKN22a}, where these schemes were applied to many different
nonconservative hyperbolic systems.

While the HSWLME model can be considered a simplification of the shallow water moment models derived in \cite{KT_SWM19}, they successfully
restore global hyperbolicity and allow for the analytical computation of steady-states. The application of the HSWLME model in more
simulations, however, is so far limited due to the lack of appropriate numerical schemes that preserve the structural properties of the
equations, for example, the steady states in the presence of friction.

In this paper, our goal is to design a second-order flux globalization based WB PCCU scheme for the HSWLME \eref{1.1}--\eref{1.4}. To this
end, we first discuss the steady-state solutions of the studied system in \S\ref{sec2}, and then, in \S\ref{sec3}, we develop the new
scheme. In \S\ref{sec5}, we demonstrate the performance of the proposed scheme on a number of numerical examples. The paper ends with
concluding remarks.

\section{Steady States}\label{sec2}
In this section, we study the steady states of the HSWLME \eref{1.1}--\eref{1.4} that satisfy \eref{1.5}, which can be equivalently
written as
\begin{equation*}
\left\{\begin{aligned}
(hu)_x&=0,\\
\Big(hu^2+\hf gh^2+\frac{1}{3}h\alpha_1^2+\ldots+\frac{1}{2N+1}h\alpha_N^2\Big)_x&=-ghZ_x+P_0,\\
(2hu\alpha_1)_x&=u(h\alpha_1)_x+P_1,\\
&\qquad\qquad~\,\vdots\\
(2hu\alpha_N)_x&=u(h\alpha_N)_x+P_{N}.
\end{aligned}\right.
\end{equation*}
The simplest ``lake-at-rest'' steady states are
\begin{equation}
\begin{aligned}
u\equiv0,\quad h+Z\equiv{\rm Const},\quad\alpha_i\equiv0,\quad i=1,\ldots,N.
\end{aligned}
\label{2.1}
\end{equation}
The non-trivial, moving-water equilibria are
\begin{equation}
\begin{aligned}
&q:=hu\equiv{\rm Const}\ne0,\quad E:=\frac{u^2}{2}+g(h+Z)+\frac{3}{2}\sum\limits_{i=1}^N\frac{\alpha_i^2}{2i+1}-{\cal P}\equiv{\rm Const},\\
&E_i:=\frac{\alpha_i}{h}-{\cal P}_i\equiv{\rm Const},\quad i=1,\ldots N,
\end{aligned}
\label{2.2}
\end{equation}
where
\begin{equation}
{\cal P}=\int\limits_{\widehat x}^x\frac{P_0(\xi,t)}{h(\xi,t)}\,{\rm d}\xi+
\sum\limits_{i=1}^N\int\limits_{\widehat x}^x\frac{P_i(\xi,t)\alpha_i(\xi,t)}{(2i+1)q(\xi,t)}\,{\rm d}\xi,\quad
{\cal P}_i=\int\limits_{\widehat x}^x\frac{P_i(\xi,t)}{h(\xi,t)q(\xi,t)}\,{\rm d}\xi,\quad i=1,\ldots,N.
\label{2.3}
\end{equation}
Notice that \eref{1.5} can be rewritten in the following form:
\begin{equation}
\bm F(\bm U)_x- Q(\bm U)\bm U_x-\bm S(\bm U)-\bm P(\bm U)=M(\bm U)\bm E(\bm U)_x=\bm0,
\label{2.4}
\end{equation}
where $\bm E:=(q,E,E_1,\ldots E_N)^\top$ are the equilibrium variables and $M$ is an $(N+2)\times(N+2)$ matrix
\begin{equation}
M(\bm U):=\begin{pmatrix}1&0&0&\ldots&0\\u&h&-\frac{1}{3}h^2\alpha_1&\ldots&-\frac{1}{2N+1}h^2\alpha_N\\2\alpha_1&0&h^2u&&\\
\vdots&&&\ddots&\\2\alpha_N&&&&h^2u\end{pmatrix}.
\label{2.5}
\end{equation}
As in \cite{KLX_21}, \eref{2.4}--\eref{2.5} are going to be used to develop a flux globalization WB PCCU scheme.

\section{Flux Globalization Based WB PCCU Scheme}\label{sec3}
In this section, we develop a second-order flux globalization based  WB PCCU scheme for the HSWLME \eref{1.1}--\eref{1.4}.

We first introduce the finite-volume cells $C_k:=\big[x_\kmh,x_\kph\big]$ of uniform (for simplicity of presentation) size
$x_\kph-x_\kmh\equiv\dx$ centered at $x_k=\big(x_\kmh+x_\kph\big)/2$, $k=1,\ldots,N_x$. Assuming that the approximate solution, realized in
terms of its cell averages $\,\xbar{\bm U}_k(t)\approx\frac{1}{\dx}\int_{C_k}\bm U(x,t)\,{\rm d}x$, is available at a certain time level
$t$, we evolve it in time by solving the following system of ODEs:
\begin{equation}
\frac{\rm d}{{\rm d}t}\,\xbar{\bm U}_k=-\frac{\bm{{\cal K}}_\kph-\bm{{\cal K}}_\kmh}{\dx},
\label{3.1}
\end{equation}
where $\bm{{\cal K}}_\kph$ are the WB PCCU fluxes, given by (see \cite{KLX_21})
\begin{equation}
\bm{{\cal K}}_\kph=\frac{a^+_\kph\bm K^-_\kph-a^-_\kph\bm K^+_\kph}{a^+_\kph-a^-_\kph}+
\frac{a^+_\kph a^-_\kph}{a^+_\kph-a^-_\kph}\Big(\widehat{\bm U}^+_\kph-\widehat{\bm U}^-_\kph\Big).
\label{3.2}
\end{equation}
Here, $\bm K^\pm_\kph=\bm K\bm(\bm U^\pm_\kph\bm)$, $\bm U^\pm_\kph$ and $\widehat{\bm U}^\pm_\kph$ are two second-order approximations of
the point values $\bm U(x_\kph,t)$, and $a_\kph^\pm$ are the one-sided local speeds of propagation, which can be estimated using the
eigenvalues of the matrix $\frac{\partial\bm F}{\partial\bm U}-Q$. Based on the eigenvalues computed in \cite{KP_SWM}, the simplest estimate
is
\begin{equation}
\hspace*{-0.2cm}\begin{aligned}
&a_\kph^+=\max\left\{u_\kph^++\sqrt{gh_\kph^++\sum\limits_{i=1}^N\frac{3\big((\alpha_i)_\kph^+\big)^2}{2i+1}},\,
u_\kph^-+\sqrt{gh_\kph^-+\sum\limits_{i=1}^N\frac{3\big((\alpha_i)_\jph^-\big)^2}{2i+1}},\,0\right\},\\
&a_\kph^-=\min\left\{u_\kph^+-\sqrt{gh_\kph^++\sum\limits_{i=1}^N\frac{3\big(\alpha_i)_\kph^+\big)^2}{2i+1}},\,
u_\kph^--\sqrt{gh_\kph^-+\sum\limits_{i=1}^N\frac{3\big((\alpha_i)_\jph^-\big)^2}{2i+1}},\,0\right\}.
\label{3.3}
\end{aligned}
\end{equation}
Notice that all of the indexed quantities in \eref{3.1}--\eref{3.3} are time-dependent, but from here on we omit this dependence for the
sake of brevity.

\subsection{Well-Balanced Reconstruction of the One-Sided Point Values}
In order to use \eref{3.2} and \eref{3.3}, one needs to compute the one-sided point values $h_\kph^\pm$, $q_\kph^\pm$,
$(\alpha_i)_\kph^\pm$, and $Z_\kph^\pm$. To make the resulting scheme WB, we need to perform piecewise linear reconstruction for the
equilibrium variables $\bm E$. To this end, we start by computing the cell center point values of $\bm E$ out of the available cell averages
$\xbar h_k$, $\xbar q_k$, $(\xbar{h\alpha_i})_k$, and the point values $Z_k:=Z(x_k)$.

First, we use \eref{2.2} and \eref{2.3} to obtain
\begin{equation*}
E_k=\frac{(u_k)^2}{2}+g(\xbar h_k+Z_k)+\frac{3}{2}\sum\limits_{i=1}^N\frac{(\alpha_i)_k^2}{2i+1}-{\cal P}_k,\quad
(E_i)_k=\frac{(\alpha_i)_k}{\xbar h_k}-({\cal P}_i)_k,\quad i=1,\ldots,N,
\end{equation*}
where $u_k=\,\xbar q_k/\,\xbar h_k$, $(\alpha_i)_k=(\xbar{h\alpha_i})_k/\,\xbar h_k$, and
\begin{equation*}
{\cal P}_k=\int\limits_{\widehat x}^{x_k}\frac{P_0(\xi,t)}{h(\xi,t)}\,{\rm d}\xi+
\sum\limits_{i=1}^N\int\limits_{\widehat x}^{x_j}\frac{P_i(\xi,t)\alpha_i(\xi,t)}{(2i+1)q(\xi,t)}\,{\rm d}\xi,\quad
({\cal P}_i)_k=\int\limits_{\widehat x}^{x_k}\frac{P_i(\xi,t)}{h(\xi,t)q(\xi,t)}\,{\rm d}\xi,\quad i=1,\ldots,N.
\end{equation*}
Notice that ${\cal P}_k$ and $({\cal P}_i)_k$ are global integral terms, we therefore use a recursive way and apply the trapezoidal rule to
evaluate them for each $[x_{k-1},x_k]$ interval for $k=2,\ldots,N_x$:
\begin{equation}
\begin{aligned}
{\cal P}_k&={\cal P}_{k-1}+\int\limits^{x_k}_{x_{k-1}}\frac{P_0(\xi,t)}{h(\xi,t)}\,{\rm d}\xi+
\sum\limits_{i=1}^N\int\limits^{x_k}_{x_{k-1}}\frac{P_i(\xi,t)\alpha_i(\xi,t)}{(2i+1)q(\xi,t)}\,{\rm d}\xi\\
&\approx{\cal P}_{k-1}+\frac{\dx}{2}\bigg[\frac{(P_0)_k}{\xbar h_k}+\frac{(P_0)_{k-1}}{\xbar h_{k-1}}\bigg]+
\frac{\dx}{2}\sum\limits_{i=1}^N\bigg[\frac{(P_i)_k(\alpha_i)_k}{{(2i+1)\,\xbar q_k}}+
\frac{(P_i)_{k-1}(\alpha_i)_{k-1}}{{(2i+1)\,\xbar q_{k-1}}}\bigg],\\
({\cal P}_i)_k&=({\cal P}_i)_{k-1}+\int\limits_{x_{k-1}}^{x_k}\frac{P_i(\xi,t)}{h(\xi,t)q(\xi,t)}\,{\rm d}\xi\approx({\cal P}_i)_{k-1}+
\frac{\dx}{2}\bigg[\frac{(P_i)_k}{\xbar h_k\,\xbar q_k}+\frac{(P_i)_{k-1}}{\xbar h_{k-1}\,\xbar q_{k-1}}\bigg],\quad i=1,\ldots,N,
\end{aligned}
\label{3.4}
\end{equation}
where
\begin{equation}
\begin{aligned}
&(P_0)_k=-\frac{\nu}{\lambda}\Big(u_k+\sum\limits_{j=1}^N(\alpha_j)_k\Big),\quad k=1,\ldots,N_x,\\
&(P_i)_k=-\frac{\nu}{\lambda}(2i+1)\Big(u_k+\sum\limits_{j=1}^N(\alpha_j)_k\Big)-
4\frac{\nu}{h_k}(2i+1)\sum\limits_{j=1}^Na_{i,j}(\alpha_j)_k,\quad i=1,\ldots,N.
\end{aligned}
\label{3.5}
\end{equation}
The values ${\cal P}_1$ and $({\cal P}_i)_1$, needed to start the iterations \eref{3.4}--\eref{3.5}, are obtained using the trapezoidal
rule, which yields
\begin{equation}
\begin{aligned}
{\cal P}_1 &=\int\limits^{x_1}_{x_\hf}\frac{P_0(\xi,t)}{h(\xi,t)}\,{\rm d}\xi+
\sum\limits_{i=1}^N\int\limits^{x_1}_{x_\hf}\frac{P_i(\xi,t)\alpha_i(\xi,t)}{(2i+1)q(\xi,t)}\,{\rm d}\xi\\
&\approx\frac{\dx}{4}\Bigg[\frac{(P_0)_1}{\xbar h_1}+\frac{(P_0)_\hf}{h_\hf}\Bigg]+
\frac{\dx}{4}\sum\limits_{i=1}^N\Bigg[\frac{(P_i)_1(\alpha_i)_1}{{(2i+1)\,\xbar q_1}}+\frac{(P_i)_\hf(\alpha_i)_\hf}{{(2i+1)q_\hf}}\Bigg],\\
({\cal P}_i)_1&=\int\limits_{x_\hf}^{x_1}\frac{P_i(\xi,t)}{h(\xi,t)q(\xi,t)}\,{\rm d}\xi\approx
\frac{\dx}{4}\Bigg[\frac{(P_i)_1}{\xbar h_1\,\xbar q_1}+\frac{(P_i)_\hf}{h_\hf q_\hf}\Bigg],\quad i=1,\ldots,N,
\end{aligned}
\label{3.6}
\end{equation}
where
\begin{equation*}
\begin{aligned}
&({\cal P}_0)_\hf=-\frac{\nu}{\lambda}\Big(u_\hf+\sum\limits_{j=1}^N(\alpha_j)_\hf\Big),\\
&({\cal P}_i)_\hf=-\frac{\nu}{\lambda}(2i+1)\Big(u_\hf+\sum\limits_{j=1}^N(\alpha_j)_\hf\Big)-
4\frac{\nu}{h_\hf}(2i+1)\sum\limits_{j=1}^Na_{i,j}(\alpha_j)_\hf,\quad i=1,\ldots,N,
\end{aligned}
\end{equation*}
and $h_\hf$, $u_\hf$, and  $(\alpha_i)_\hf$, $i=1,\ldots,N$, are determined based on the prescribed boundary conditions.

Equipped with the point values $E_k$ and $(E_i)_k$, $i=1,\ldots,N$, we proceed with performing the generalized minmod reconstruction
(provided in Appendix \ref{appxA}) for the equilibrium variables $q$, $E$, and $E_i$, $i=1,\ldots,N$, as well as for $Z$ to obtain the
left/right-sided point values $q_\kph^\pm$, $E_\kph^\pm$,  $(E_i)_\kph^\pm$, $i=1,\ldots,N$, and $Z_\kph^\pm$. We then use Newton's method
to numerically solve the following nonlinear equations:
\begin{equation}
E_\kph^\pm=\frac{\big(q_\kph^\pm\big)^2}{2\big(h_\kph^\pm\big)^2}+g\big(h_\kph^\pm+Z_\kph^\pm\big)+
\frac{3}{2}\sum\limits_{i=1}^N\frac{\left(\big[(E_i)_\kph^\pm+({\cal P}_i)_\kph\big]h_\kph^\pm\right)^2}{2i+1}-{\cal P}_\kph
\label{3.7}
\end{equation}
for $h_\kph^\pm$. Once the point values $h_\kph^\pm$ are obtained, we compute $u_\kph^\pm=q_\kph^\pm/h_\kph^\pm$, and
\begin{equation}
(\alpha_i)_\kph^\pm=\big[(E_i)_\kph^\pm+({\cal P}_i)_\kph\big]h_\kph^\pm,
\label{3.8f}
\end{equation}
from \eref{2.2}.

In \eref{3.7}, the point values ${\cal P}_\kph$ and $({\cal P}_i)_\kph$, $i=1,\ldots,N$, are obtained using the midpoint rule applied in the
following recursive way for $k=1,\ldots,N_x$:
\begin{equation}
\begin{aligned}
{\cal P}_\kph&={\cal P}_\kmh+\int\limits^{x_\kph}_{x_\kmh}\frac{P_0(\xi,t)}{h(\xi,t)}\,{\rm d}\xi+
\sum\limits_{i=1}^N\int\limits^{x_\kph}_{x_\kmh}\frac{P_i(\xi,t)\alpha_i(\xi,t)}{(2i+1)q(\xi,t)}\,{\rm d}\xi\\
&\approx{\cal P}_\kmh+\frac{(P_0)_k}{\xbar h_k}\,\dx+\sum\limits_{i=1}^N\frac{(P_i)_k(\alpha_i)_k}{{(2i+1)\,\xbar q_k}}\,\dx,\\
({\cal P}_i)_\kph&=({\cal P}_i)_\kmh+\int\limits_{x_{k-1}}^{x_k}\frac{P_i(\xi,t)}{h(\xi,t)q(\xi,t)}\,{\rm d}\xi\approx
({\cal P}_i)_\kmh+\frac{(P_i)_k}{\xbar h_k\,\xbar q_k}\,\dx,\quad i=1,\ldots, N,
\end{aligned}
\label{3.8}
\end{equation}
with ${\cal P}_\hf=({\cal P}_i)_\hf=0$.
\begin{rmk}
In \eref{3.4}, \eref{3.6}, and \eref{3.8}, one has to divide by $\xbar h_k$ and $\xbar q_k$. These divisions have to be desingularized. The
desingularization is performed as follows. If we need to divide, say, $B$ by $A$, we replace $\frac{B}{A}$ with
$\frac{2AB}{A^2+[\max(|A|,\varepsilon)]^2}$, where $\varepsilon$ is a small positive number taken to be $\varepsilon=10^{-6}$ in the
numerical experiments reported in \S\ref{sec5}. This ensures stability of the computations in local near-dry situations of small $h$ and/or small discharge $q$.
\end{rmk}

\subsection{Well-Balanced Global Flux Evaluation}
We now proceed with the evaluation of the one-sided values of the global fluxes given by
\begin{equation*}
\bm K^\pm_\kph=\bm F\big(\bm U^\pm_\kph\big)-\bm R^\pm_\kph,
\end{equation*}
which, according to \eref{1.2} and \eref{1.7}, are
\begin{equation*}
\bm K^\pm_\kph=\begin{pmatrix}q^\pm_\kph\\[1.2ex]K^\pm_\kph\\[1.5ex](K^{(1)})^\pm_\kph\\\vdots\\(K^{(N)})^\pm_\kph\end{pmatrix}=
\begin{pmatrix}q^\pm_\kph\\[1.2ex]q^\pm_\kph u^\pm_\kph+\dfrac{g}{2}\big(h^\pm_\kph\big)^2+
\sum\limits_{i=1}^N\dfrac{1}{2i+1}h^\pm_\kph\big((\alpha_i)^\pm_\kph\big)^2-R^\pm_\kph\\[1.5ex]
2q^\pm_\kph(\alpha_i)^\pm_\kph-(R^{(1)})^\pm_\kph\\\vdots\\2q^\pm_\kph(\alpha_N)^\pm_\kph-(R^{(N)})^\pm_\kph\end{pmatrix},
\end{equation*}
where
\begin{equation*}
\begin{aligned}
&R=-\int\limits_{\widehat x}^x\big[gh(\xi,t)Z_\xi(\xi,t)-P_0(\xi,t)\big]{\rm d}\xi,\\
&R^{(i)}=\int\limits_{\widehat x}^x\big\{u(\xi,t)[h(\xi,t)\alpha_i(\xi,t)]_\xi+P_i(\xi,t)\big\}{\rm d}\xi,\quad i=1,\ldots,N.
\end{aligned}
\end{equation*}
Notice that the global terms $R$ and $R^{(i)}$, $i=1,\ldots,N$, contain integrals of nonconservative products and we therefore use the
path-conservative integration technique introduced in \cite{KLX_21} to evaluate $R^\pm_\kph$ and $(R^{(i)})^\pm_\kph$, $i=1,\ldots,N$, in
the following recursive way:
$$
\begin{aligned}
&R_\hf^-=0,\quad R_\hf^+=B_{\bm\Psi,\hf}\approx-\int\limits_0^1gh_\hf(s)Z'_\hf(s)\,{\rm d}s,\\
&R^-_\kph=R^+_\kmh+B_k\approx R^+_\kmh-\int\limits_{C_k}(ghZ_x-P_0)\,{\rm d}x,\quad k=1,\ldots,N_x,\\
&R^+_\kph=R^-_\kph+B_{\bm\Psi,\kph}\approx R^-_\kph-\int\limits^1_0gh(s)Z'(s)\,{\rm d}s,\quad k=1,\ldots,N_x.
\end{aligned}
$$
and
$$
\begin{aligned}
&(R^{(i)})_\hf^-=0,\quad(R^{(i)})_\hf^+=B_{\bm\Psi,\hf}^{(i)}\approx\int\limits_0^1u(s)[h(s)\alpha_i(s)]'\,{\rm d}s,\\
&(R^{(i)})^-_\kph=(R^{(i)})^+_\kmh+B_k^{(i)}\approx(R^{(i)})^+_\kmh-\int\limits_{C_k}[u(h\alpha_i)_x+P_i]\,{\rm d}x,\quad k=1,\ldots,N_x,\\
&(R^{(i)})^+_\kph=(R^{(i)})^-_\kph+B_{\bm\Psi,\kph}^{(i)}\approx(R^{(i)})^-_\kph-\int\limits^1_0u(s)[h(s)\alpha_i(s)]'\,{\rm d}s,\quad
k=1,\ldots,N_x,
\end{aligned}
$$
for $i=1,\ldots,N$.
Here, $\bm\Psi_\kph(s):=\bm\Psi\big(s;\bm U^-_\kph,\bm U^+_\kph\big)$ is a sufficiently smooth path connecting the states $\bm U^-_\kph$ and
$\bm U^+_\kph$, that is,
\begin{equation*}
\bm\Psi:[0,1]\times\mathbb R^N\times\mathbb R^N\to\mathbb R^N,\quad\bm\Psi\big(0;\bm U^-_\kph,\bm U^+_\kph\big)=\bm U^-_\kph,\quad
\bm\Psi\big(1;\bm U^-_\kph,\bm U^+_\kph\big)=\bm U^+_\kph,
\end{equation*}
where $\bm\Psi_\kph(s)$ should be obtained in a WB manner by connecting the left and right states of the equilibrium variables using a
linear segment path
\begin{equation*}
\bm E_\kph(s)=\bm E^-_\kph+s\big(\bm E^+_\kph-\bm E^-_\kph\big),\quad s\in[0,1].
\end{equation*}
Following the technique introduced in \cite{KLX_21,CKLX_22}, we design WB trapezoidal-like quadratures for
$\bm B_k=\big(0,B_k,B_k^{(1)},\ldots,B_k^{(N)}\big)^\top$ and
$\bm B_{\bm\Psi,\kph}=\big(0,B_{\bm\Psi,\kph},B_{\bm\Psi,\kph}^{(1)},\ldots,B_{\bm\Psi,\kph}^{(N)}\big)^\top$ based on \eref{2.4} and
\eref{2.5}:
\begin{equation*}
\begin{aligned}
\bm B_k&=\bm F\big(\bm U_\kph^-\big)-\bm F\big(\bm U_\kmh^+\big)-\hf\left[M\big(\bm U_\kph^-\big)+M\big(\bm U_\kmh^+\big)\right]
\left(\bm E_\kph^--\bm E_\kmh^+\right),\\
\bm B_{\bm\Psi,\kph}&=\bm F\big(\bm U_\kph^+\big)-\bm F\big(\bm U_\kph^-\big)-
\hf\left[M\big(\bm U_\kph^+\big)+M\big(\bm U_\kph^-\big)\right]\left(\bm E_\kph^+-\bm E_\kph^-\right).
\end{aligned}
\end{equation*}
The nonzero components of $\bm B_k$ and $\bm B_{\bm\Psi,\kph}$ read as
\allowdisplaybreaks
\begin{align*}
B_k&=q_\kph^-u_\kph^--q_\kmh^+u_\kmh^++\frac{g}{2}\Big[\big(h^-_\kph\big)^2-\big(h^+_\kmh\big)^2\Big]\\
&+\sum\limits_{i=1}^N\frac{1}{2i+1}\left[h^-_\kph\big((\alpha_i)^-_\kph\big)^2-h^+_\kmh\big((\alpha_i)^+_\kmh\big)^2\right]\\
&-\hf\big(h_\kph^-+h_\kmh^+\big)\big(E_\kph^--E_\kmh^+\big)-\hf\big(u_\kph^-+u_\kmh^+\big)\big(q_\kph^--q_\kmh^+\big)\\
&+\hf\sum\limits_{i=1}^N\frac{1}{2i+1}\left[\big(h^-_\kph\big)^2(\alpha_i)^-_\kph+\big(h^+_\kmh\big)^2(\alpha_i)^+_\kmh\right]
\big((E_i)^-_\kph-(E_i)^+_\kmh\big),\\
B_{\bm\Psi,\kph}&=q_\kph^+u_\kph^+-q_\kph^-u_\kph^-+\frac{g}{2}\Big[\big(h^+_\kph\big)^2-\big(h^-_\kph\big)^2\Big]\\
&+\sum\limits_{i=1}^N\frac{1}{2i+1}\left[h^+_\kph\big((\alpha_i)^+_\kph\big)^2-h^-_\kph\big((\alpha_i)^-_\kph\big)^2\right]\\
&-\hf\big(h_\kph^++h_\kph^-\big)\big(E_\kph^+-E_\kph^-\big)-\hf\big(u_\kph^++u_\kph^-\big)\big(q_\kph^+-q_\kph^-\big)\\
&+\hf\sum\limits_{i=1}^N\frac{1}{2i+1}\left[\big(h^+_\kph\big)^2(\alpha_i)^+_\kph+\big(h^-_\kph\big)^2(\alpha_i)^-_\kph\right]
\big((E_i)^+_\kph-(E_i)^-_\kph\big),
\end{align*}
and
\begin{equation*}
\begin{aligned}
B_k^{(i)}&=2q_\kph^-(\alpha_i)_\kph^--2q_\kmh^+(\alpha_i)_\kmh^+-\hf\big((\alpha_i)_\kph^-+(\alpha_i)_\kmh^+\big)
\big(q_\kph^--q_\kmh^+\big),\\
&-\hf\left[(h^-_\kph)^2u^-_\kph+(h^+_\kmh)^2u^+_\kmh\right]\big((E_i)^-_\kph-(E_i)^+_\kmh\big),\\
B_{\bm\Psi,\kph}^{(i)}&=2q_\kph^+(\alpha_i)_\kph^+-2q_\kph^-(\alpha_i)_\kph^--\hf\big((\alpha_i)_\kph^++(\alpha_i)_\kph^-\big)
\big(q_\kph^+-q_\kph^-\big),\\
&-\hf\left[(h^+_\kph)^2u^+_\kph+(h^-_\kmh)^2u^-_\kmh\right]\big((E_i)^+_\kph-(E_i)^-_\kmh\big),\quad i=1,\ldots,N.
\end{aligned}
\end{equation*}
\begin{rmk}
One can show that at steady states satisfying either \eref{2.1} or \eref{2.2}, $\bm K_\kmh^+=\bm K_\kph^-=\bm K_\kph^+$ for all $k$. This
immediately follows from \cite[Theorem 4.1]{KLX_21}.
\end{rmk}

\subsection{Well-Balanced Evolution}
Equipped with the reconstructed one-sided point values $h_\kph^\pm$, $q_\kph^\pm$, and $\bm K_\kph^\pm$, the semi-discrete flux
globalization based WB PCCU scheme \eref{3.1}--\eref{3.2} for the HSWLME \eref{1.1}--\eref{1.4} can be written as \eref{3.1} with
${\cal K}_\kph=\big(K_\kph^{(1)},\ldots,K_\kph^{(N+2)}\big)^\top$, where
\begin{equation}
\begin{aligned}
{\cal K}^{(1)}_\kph&=\frac{a^+_\kph q^-_\kph-a^-_\kph q^+_\kph}{a^+_\kph-a^-_\kph}+
\frac{a^+_\kph a^-_\kph}{a^+_\kph-a^-_\kph}\left(\widehat h^+_\kph-\widehat h^-_\kph\right),\\
{\cal K}^{(2)}_\kph&=\frac{a^+_\kph K^-_\kph-a^-_\kph K^+_\kph}{a^+_\kph-a^-_\kph}+
\frac{a^+_\kph a^-_\kph}{a^+_\kph-a^-_\kph}\left(\widehat q_\kph^{\,+}-\widehat q^{\,-}_\kph\right),\\
{\cal K}^{(i+2)}_\kph&=\frac{a^+_\kph(K^{(i)})^-_\kph-a^-_\kph(K^{(i)})^+_\kph}{a^+_\kph-a^-_\kph}\\
&+\frac{a^+_\kph a^-_\kph}{a^+_\kph-a^-_\kph}
\left(\widehat h^+_\kph(\widehat{\alpha_i})^+_\kph-\widehat h^-_\kph(\widehat{\alpha_i})^-_\kph\right),\quad i=1,\ldots,N.
\end{aligned}
\label{3.10}
\end{equation}
Here, $\widehat q^{\,\pm}_\kph=q_\kph^\pm$ as $q_\kph^+=q_\kph^-$ at the steady states \eref{2.1} and \eref{2.2}, while
$\widehat h^+_\kph$ and $\widehat h_\kph^-$ are obtained by numerically solving the following nonlinear equations (compare them with
\eref{3.7}):
\begin{equation*}
E_\kph^\pm=\frac{\big(q_\kph^\pm\big)^2}{2\big(\widehat h_\kph^\pm\big)^2}+g\big(\widehat h^\pm_\kph+Z_\kph\big)+
\frac{3}{2}\sum\limits_{i=1}^N\frac{\left(\big[(E_i)_\kph^\pm+({\cal P}_i)_\kph\big]\widehat h_\kph^\pm\right)^2}{2i+1}-{\cal P}_\kph,
\end{equation*}
where the same values $Z_\kph=\big(Z_\kph^++Z_\kph^-\big)/2$ are used in both the ``$+$'' and ``$-$'' equations. Finally,
$(\widehat{\alpha_i})^\pm_\kph$ are obtain by (compare them with \eref{3.8f})
\begin{equation*}
(\widehat{\alpha_i})_\kph^\pm=\big[(E_i)_\kph^\pm+({\cal P}_i)_\kph\big]\,\widehat h_\kph^\pm,
\end{equation*}
Notice that unlike $h_\kph^\pm$ and $(\alpha_i)^\pm_\kph$, which may be different at steady states, $\widehat h^+_\kph=\widehat h_\kph^-$
and $(\widehat{\alpha_i})^-_\kph=(\widehat{\alpha_i})^+_\kph$, $i=1,\ldots,N$, when the computed solution satisfies either \eref{2.1} or
\eref{2.2}.

\section{Numerical Examples}\label{sec5}
In this section, we demonstrate the performance of the proposed second-order semi-discrete flux globalization based WB PCCU scheme on a
number of  numerical examples. In all of the examples, we take the generalized minmod parameter $\theta=1.3$, the friction parameter
$\lambda=1$, and use nonreflecting boundary conditions. The time-dependent ODE system \eref{3.1} is integrated using the three-stage
third-order strong stability preserving (SSP) Runge-Kutta method (see, e.g., \cite{GKS,GST}) with a time step computed at every time level
using the CFL number $0.45$, namely, by taking
\begin{equation*}
\dt=0.45\,\frac{\dx}{a_{\max}},\quad a_{\max}:=\max_k\left\{\max\big(a^+_\kph,-a^-_\kph\big)\right\}.
\end{equation*}

\subsubsection*{Example 1---Convergence to Steady States}
In the first example, we study the convergence in time of the numerical solutions computed by the flux globalization based WB and
non-well-balanced (NWB) PCCU schemes towards steady flows over a hump. The proposed WB PCCU scheme uses both $h_\kph^\pm$ and
$\widehat h_\kph^\pm$, while the NWB PCCU scheme uses $h_\kph^\pm$ only, that is, it is obtained from the WB PCCU scheme by setting
$\widehat h_\kph^\pm=h_\kph^\pm$ and $(\widehat{\alpha_i})_\kph^\pm=(\widehat{\alpha_i})_\kph^\pm$, $i=1,\ldots,N$ in \eref{3.10}.

We consider both continuous and discontinuous bottom topographies given by
\begin{equation*}
Z(x)=\left\{\begin{aligned}
&0.2-0.05(x-10)^2&&\mbox{if}~8\le x\le12,\\
&0&&\mbox{otherwise},
\end{aligned}\right.
\end{equation*}
and
\begin{equation}
Z(x)=\left\{\begin{aligned}
&0.2&&\mbox{if}~8\le x\le12,\\
&0&&\mbox{otherwise},
\end{aligned}\right.
\label{4.1}
\end{equation}
respectively, and we take $N=2$ and $g=9.812$ in this example. Depending on the initial and boundary conditions, the flow may be
supercritical (Case (a)), subcritical (Case (b)), or transcritical (Case (c)). We take the following initial and
boundary conditions:
$$
\begin{aligned}
&\mbox{Case (a):}&&
\left\{\begin{aligned}
&h(x,0)=2-Z(x),&&q(x,0)=0,&&\alpha_1(x,0)=0,&&\alpha_2(x,0)=0,\\
&h(0,t)=2,&&q(0,t)=24,&&\alpha_1(0,t)=0.5,&&\alpha_2(0,t)=0.5;
\end{aligned}\right.\\
&\mbox{Case (b):}&&
\left\{\begin{aligned}
&h(x,0)=2-Z(x),&&q(x,0)=0,&&\alpha_1(x,0)=0,&&\alpha_2(x,0)=0,\\
&q(0,t)=4.42,&&h(25,t)=2,&&\alpha_1(0,t)=0.5,&&\alpha_2(0,t)=0.5;
\end{aligned}\right.\\
&\mbox{Case (c):}&&
\left\{\begin{aligned}
&h(x,0)=0.93-Z(x),&&q(x,0)=0,&&\alpha_1(x,0)=0,&&\alpha_2(x,0)=0,\\
&q(0,t)=4.42,&&h(25,t)=0.93,&&\alpha_1(0,t)=0.5,&&\alpha_2(0,t)=0.5.
\end{aligned}\right.
\end{aligned}
$$

\begin{figure}[ht!]
\centerline{\includegraphics[trim=0.0cm 0.3cm 1.cm 0.1cm,clip,width=5cm]{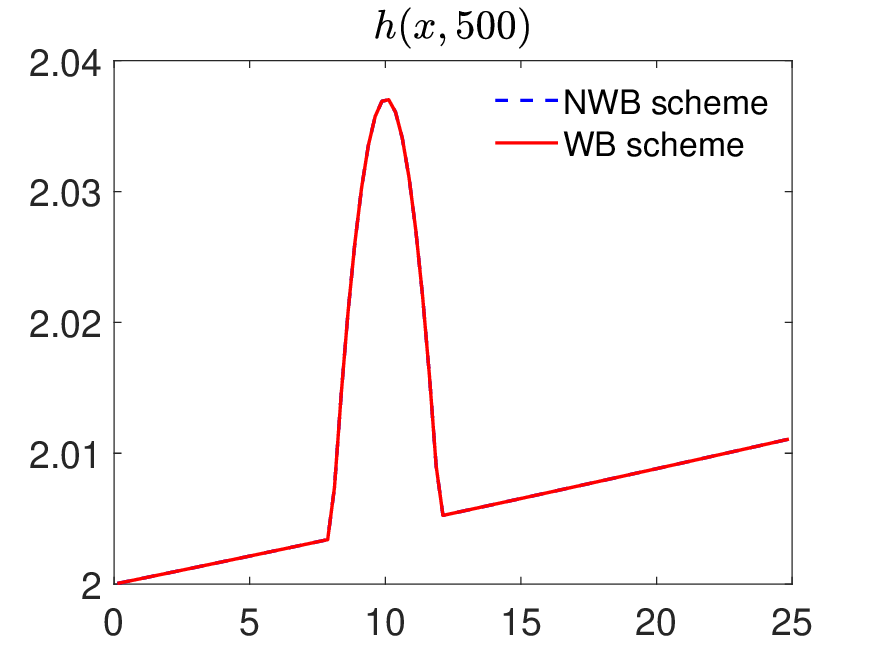}\hspace*{0.5cm}
	    \includegraphics[trim=0.0cm 0.3cm 1.cm 0.1cm,clip,width=5cm]{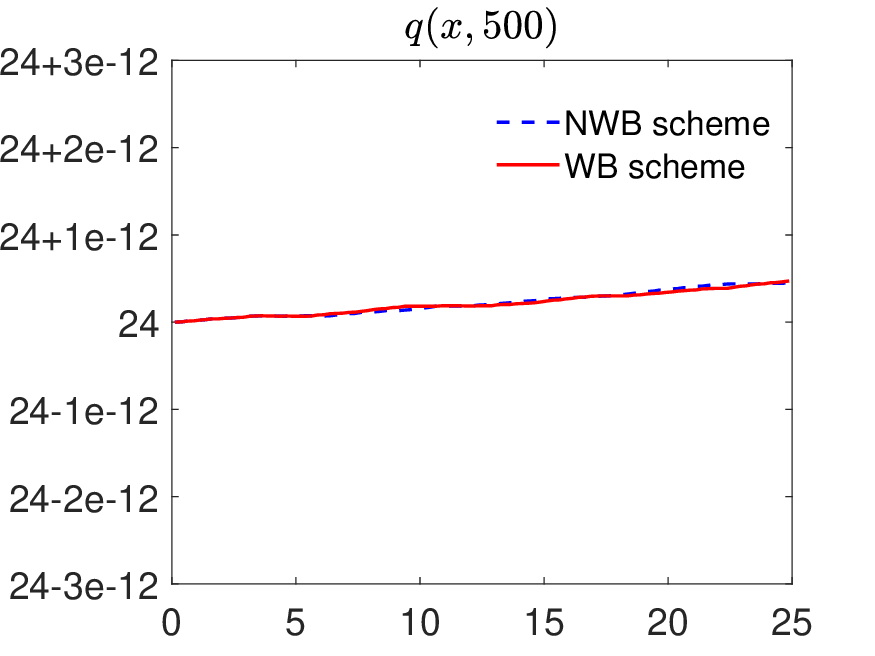}}
\vskip5pt
\centerline{\includegraphics[trim=0.0cm 0.3cm 1.cm 0.1cm,clip,width=5cm]{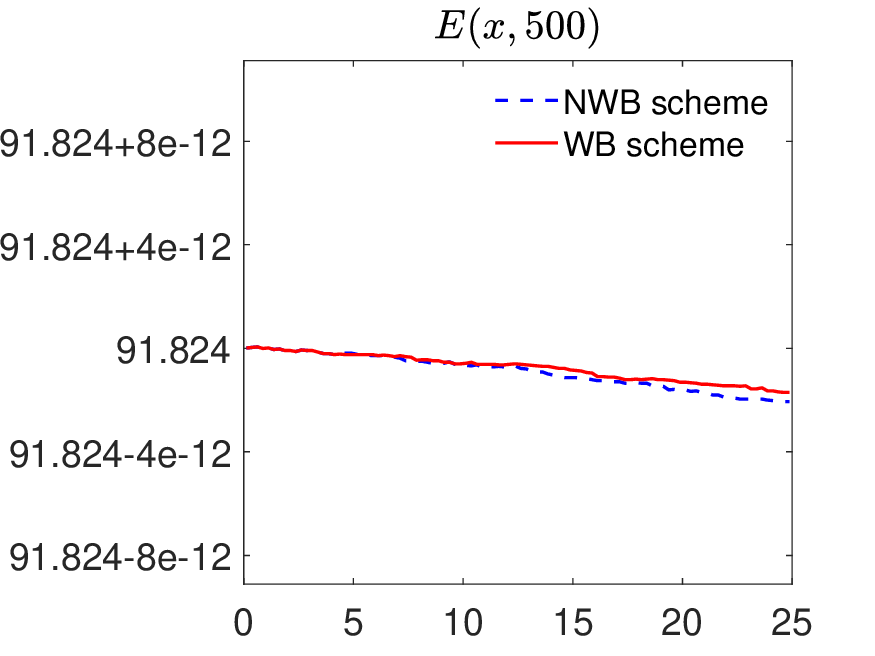}\hspace*{0.5cm}
            \includegraphics[trim=0.0cm 0.3cm 1.cm 0.1cm,clip,width=5cm]{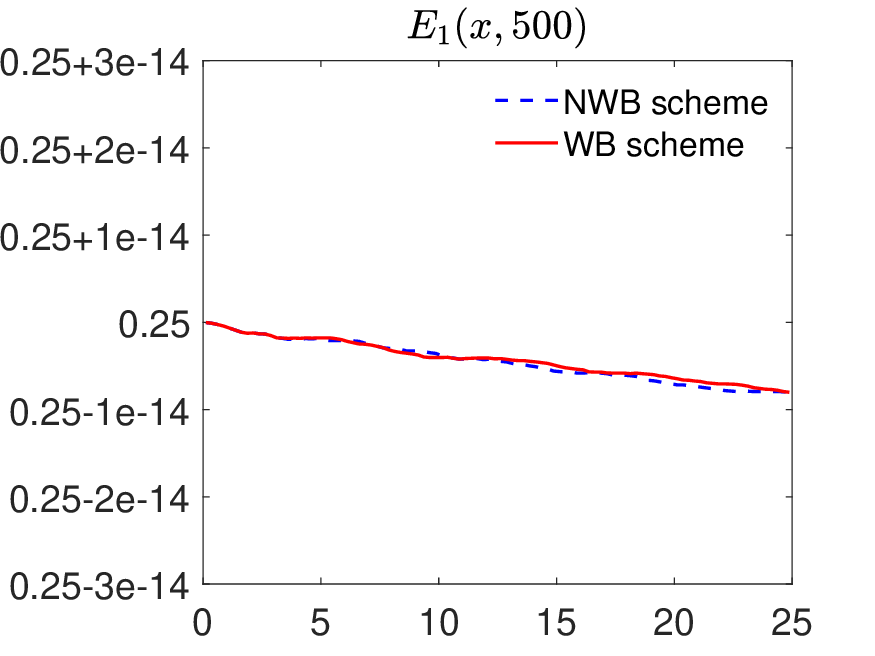}\hspace*{0.5cm}
            \includegraphics[trim=0.0cm 0.3cm 1.cm 0.1cm,clip,width=5cm]{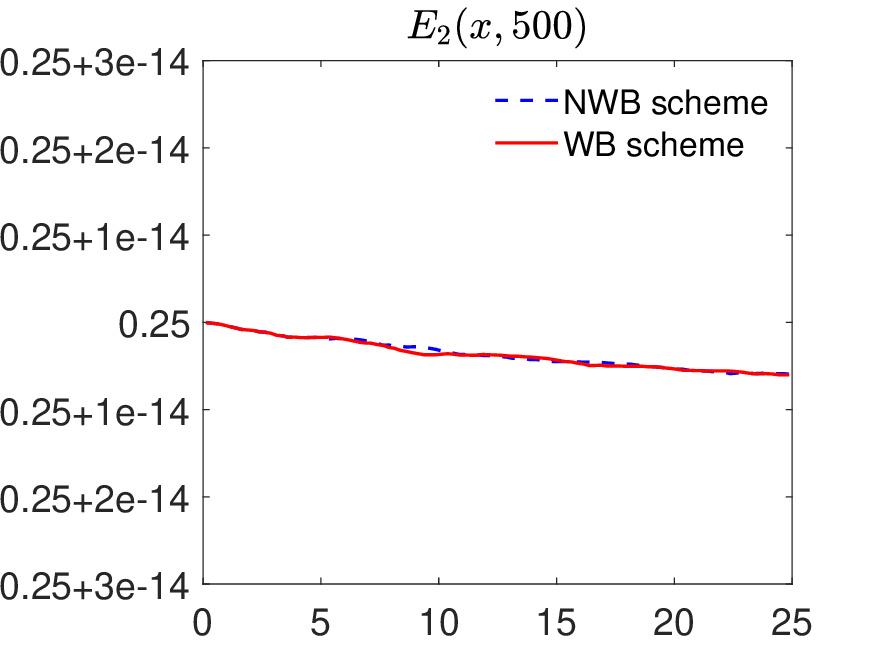}}
\vskip8pt
\centerline{\includegraphics[trim=0.0cm 0.3cm 1.cm 0.1cm,clip,width=5cm]{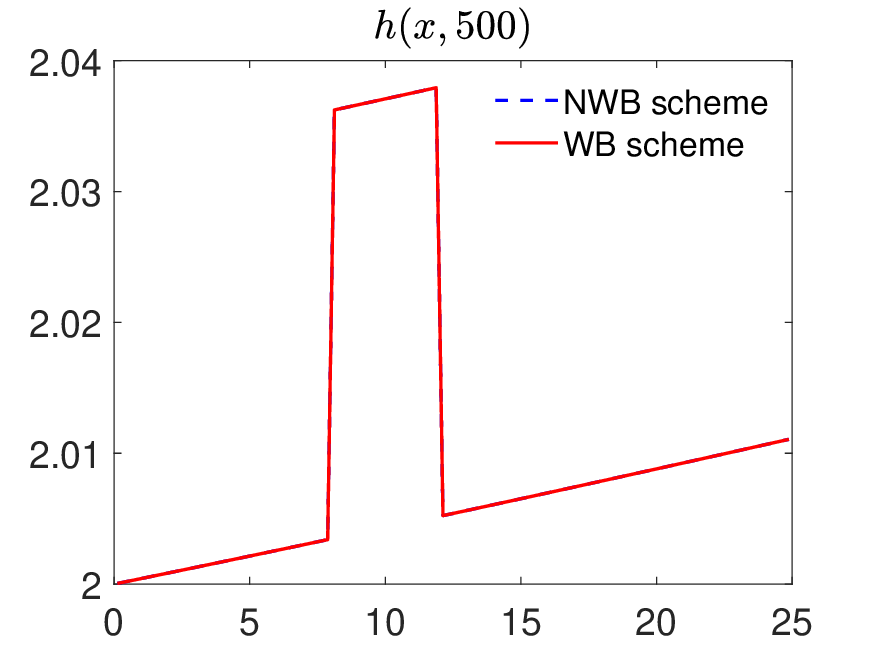}\hspace*{0.5cm}
            \includegraphics[trim=0.0cm 0.3cm 1.cm 0.1cm,clip,width=5cm]{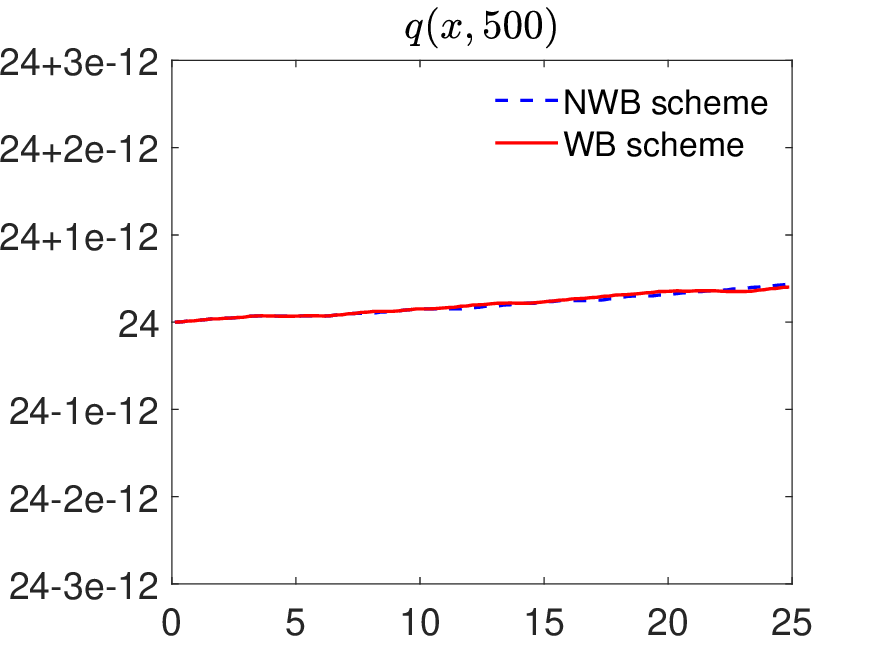}}
\vskip5pt
\centerline{\includegraphics[trim=0.0cm 0.3cm 1.cm 0.1cm,clip,width=5cm]{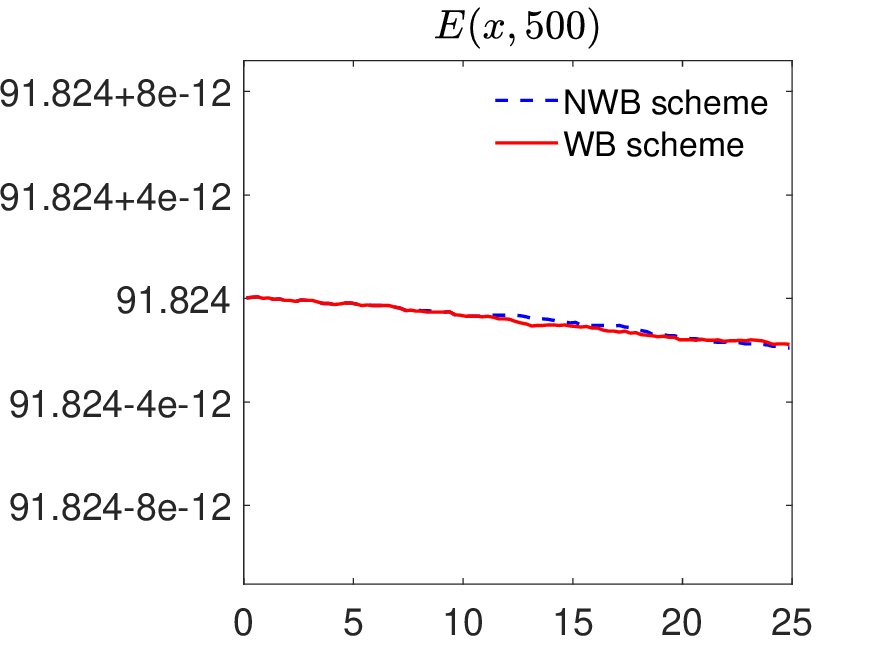}\hspace*{0.5cm}
            \includegraphics[trim=0.0cm 0.3cm 1.cm 0.1cm,clip,width=5cm]{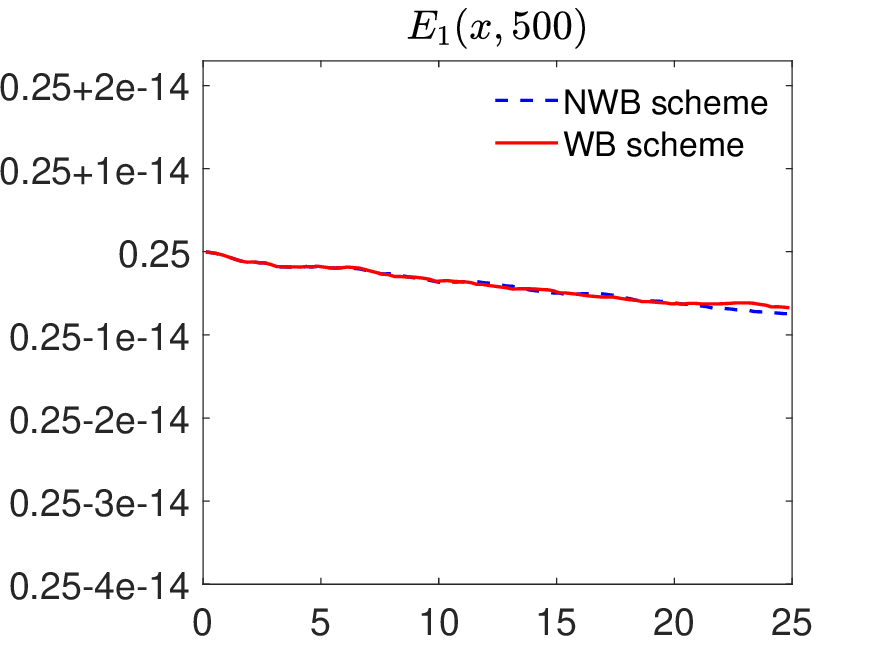}\hspace*{0.5cm}
            \includegraphics[trim=0.0cm 0.3cm 1.cm 0.1cm,clip,width=5cm]{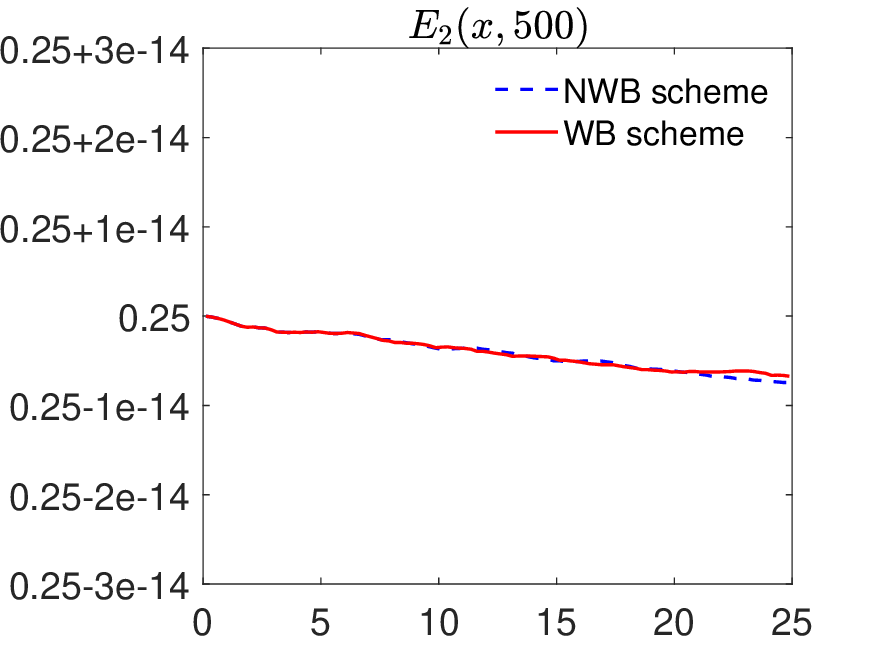}}
\caption{\sf Example 1, Case (a): $h$, $q$, $E$, $E_1$, and $E_2$ computed by the WB and NWB PCCU schemes over the continuous (top two rows)
and discontinuous (bottom two rows) $Z$.\label{fig41}}
\end{figure}

\begin{figure}[ht!]
\centerline{\includegraphics[trim=0.0cm 0.3cm 1.cm 0.1cm,clip,width=5cm]{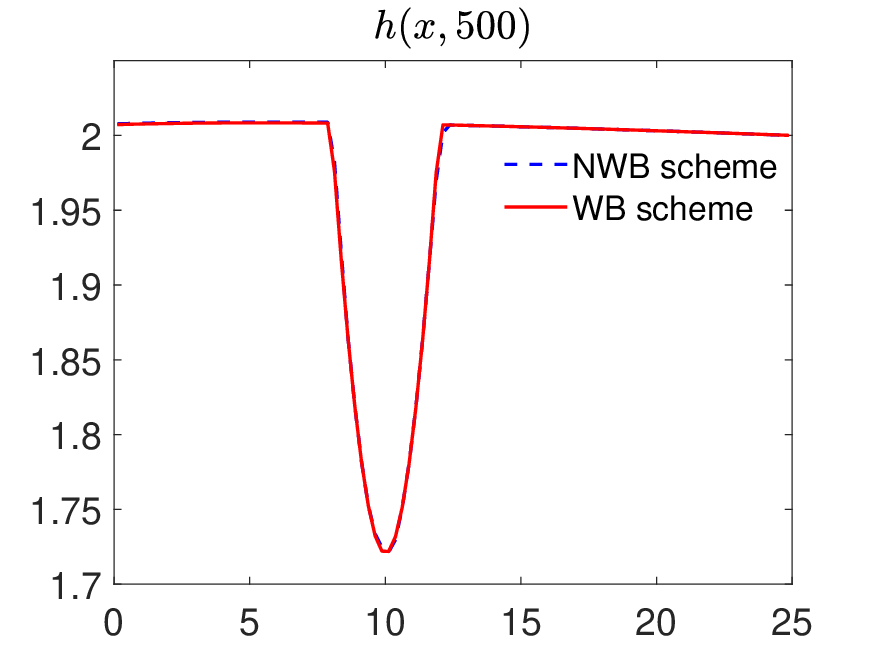}\hspace*{0.5cm}
            \includegraphics[trim=0.0cm 0.3cm 1.cm 0.1cm,clip,width=5cm]{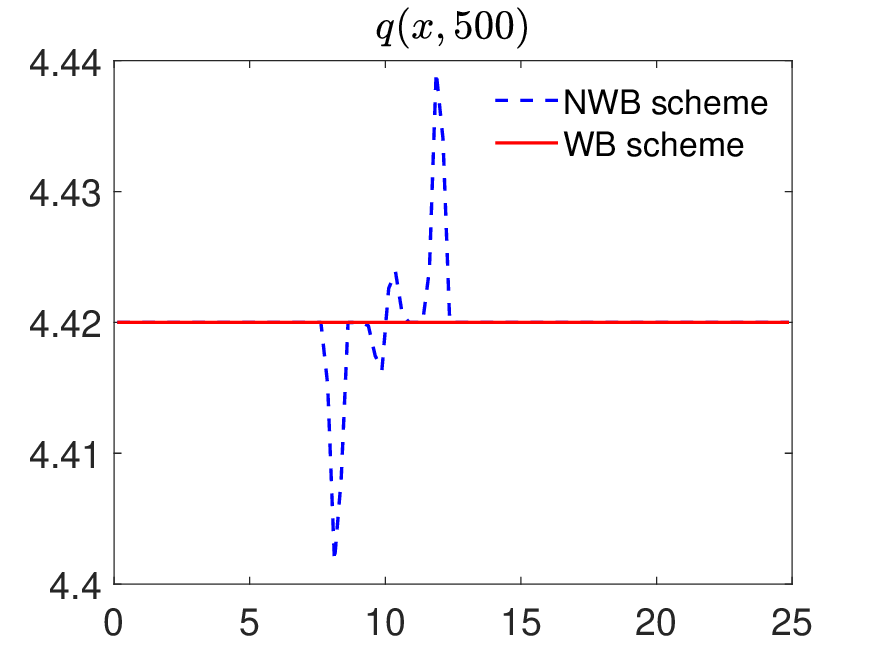}}
\vskip5pt
\centerline{\includegraphics[trim=0.0cm 0.3cm 1.cm 0.1cm,clip,width=5cm]{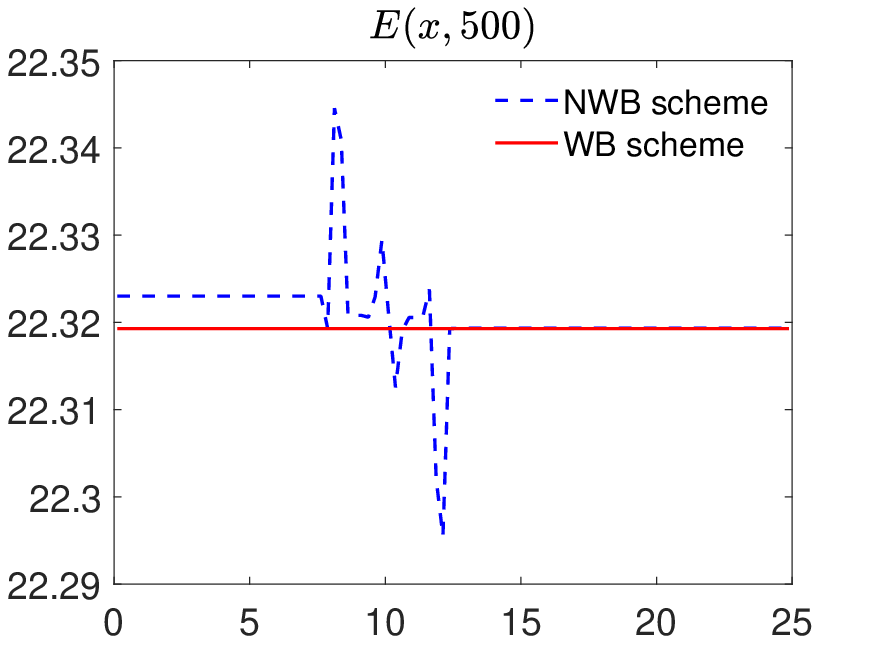}\hspace*{0.5cm}
            \includegraphics[trim=0.0cm 0.3cm 1.cm 0.1cm,clip,width=5cm]{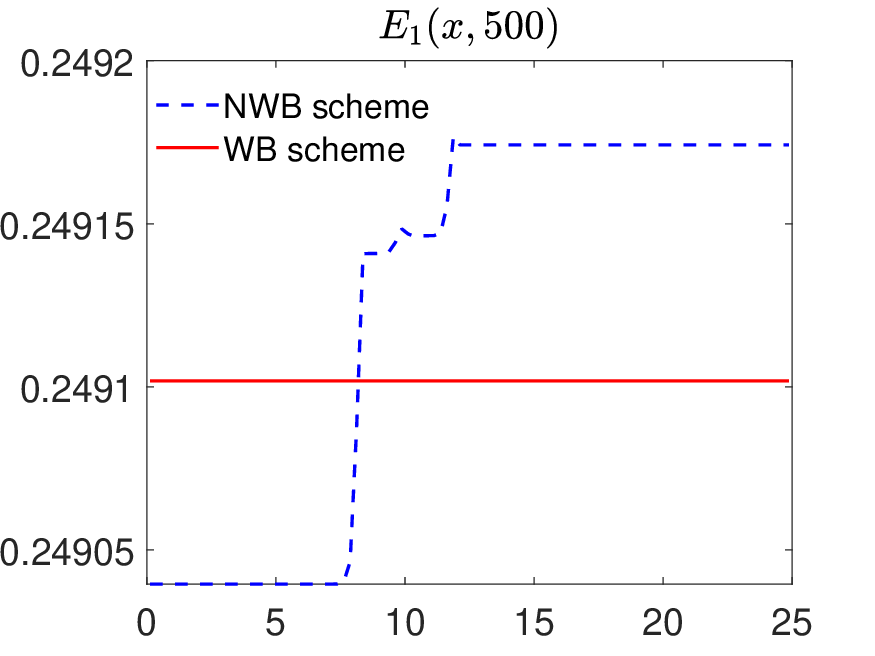}\hspace*{0.5cm}
            \includegraphics[trim=0.0cm 0.3cm 1.cm 0.1cm,clip,width=5cm]{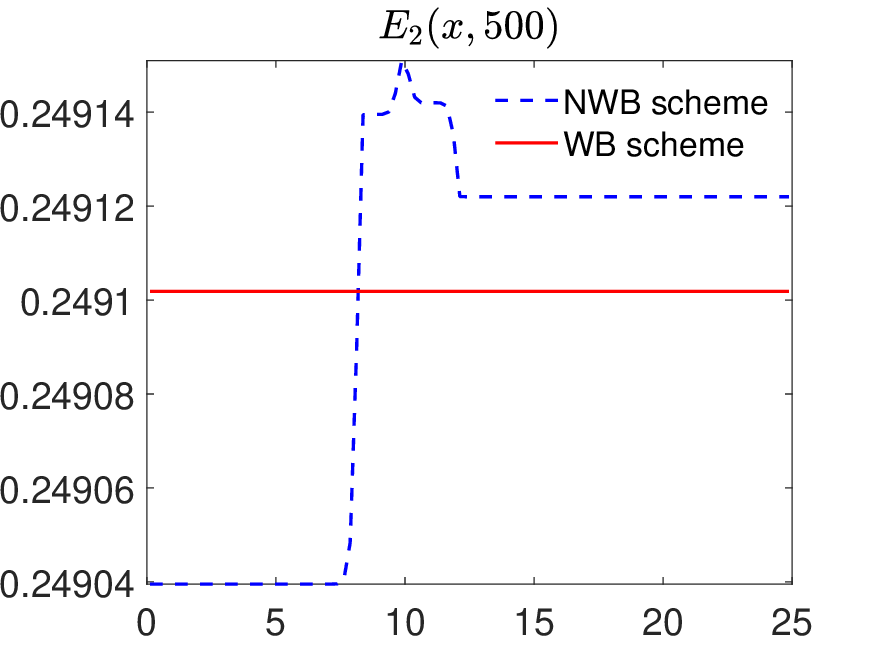}}
\vskip8pt
\centerline{\includegraphics[trim=0.0cm 0.3cm 1.cm 0.1cm,clip,width=5cm]{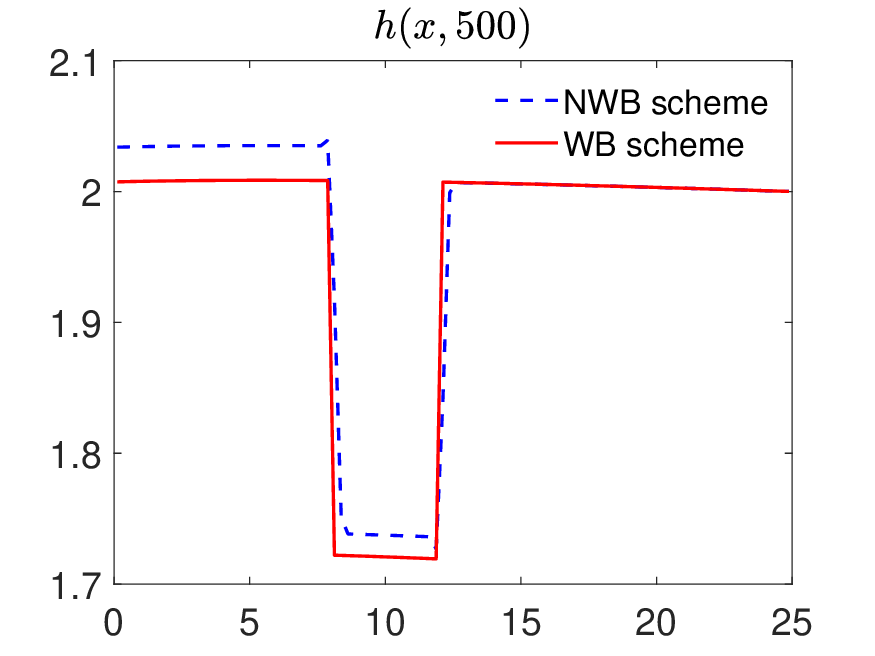}\hspace*{0.5cm}
            \includegraphics[trim=0.0cm 0.3cm 1.cm 0.1cm,clip,width=5cm]{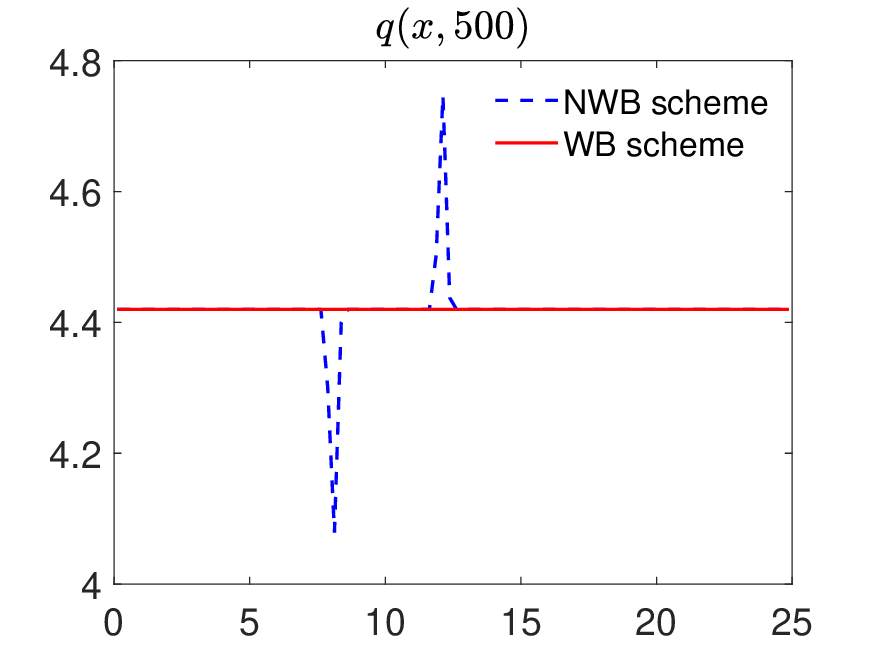}}
\vskip5pt
\centerline{\includegraphics[trim=0.0cm 0.3cm 1.cm 0.1cm,clip,width=5cm]{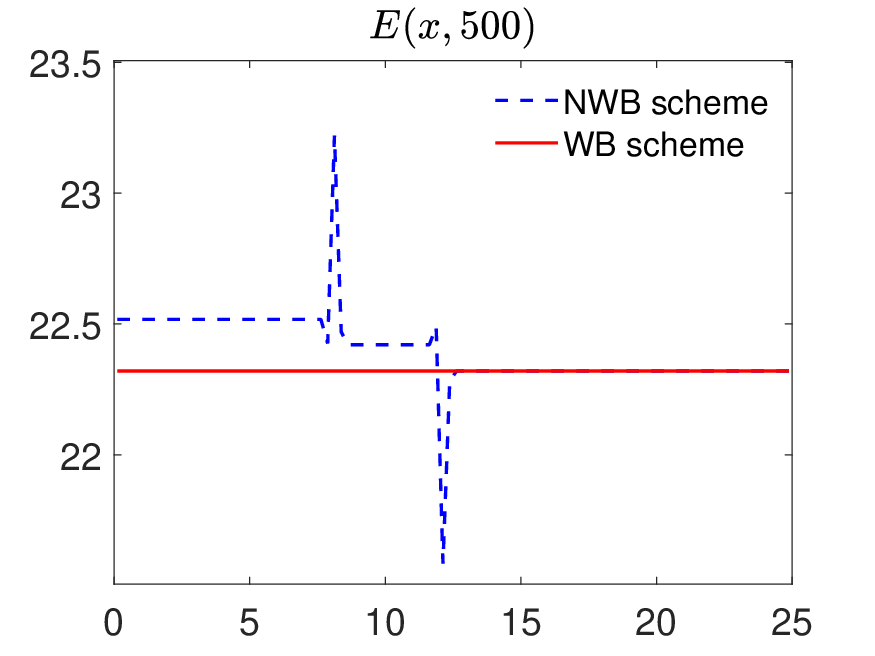}\hspace*{0.5cm}
            \includegraphics[trim=0.0cm 0.3cm 1.cm 0.1cm,clip,width=5cm]{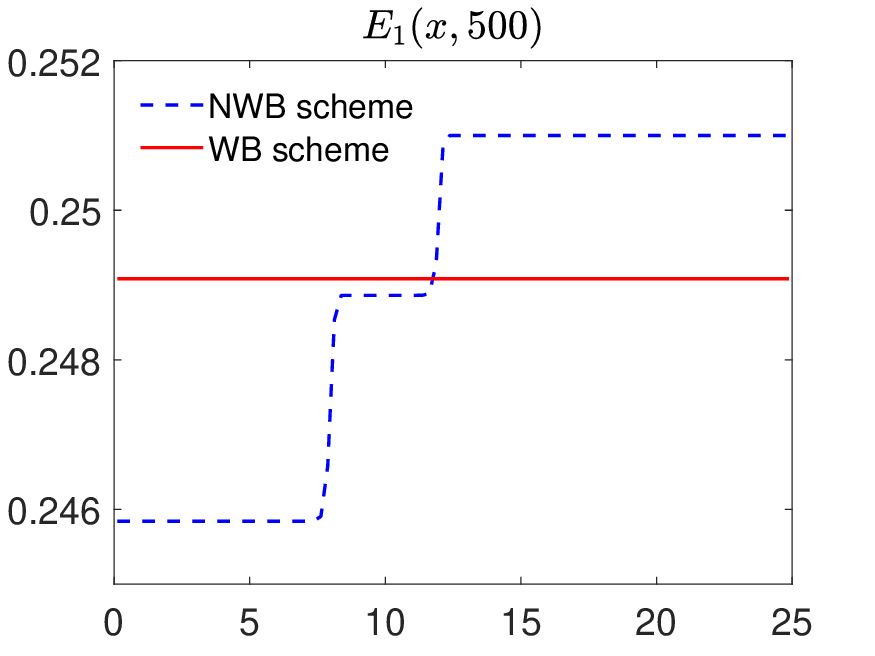}\hspace*{0.5cm}
            \includegraphics[trim=0.0cm 0.3cm 1.cm 0.1cm,clip,width=5cm]{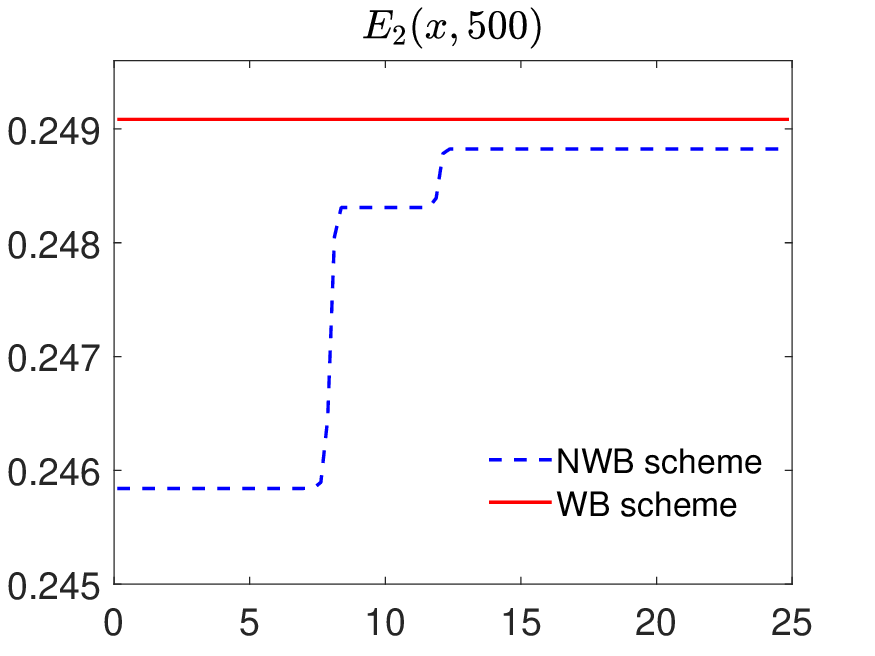}}
\caption{\sf Same as in Figure \ref{fig41}, but for Case (b).\label{fig42}}
\end{figure}

\begin{figure}[ht!]
\centerline{\includegraphics[trim=0.0cm 0.3cm 1.cm 0.1cm,clip,width=5cm]{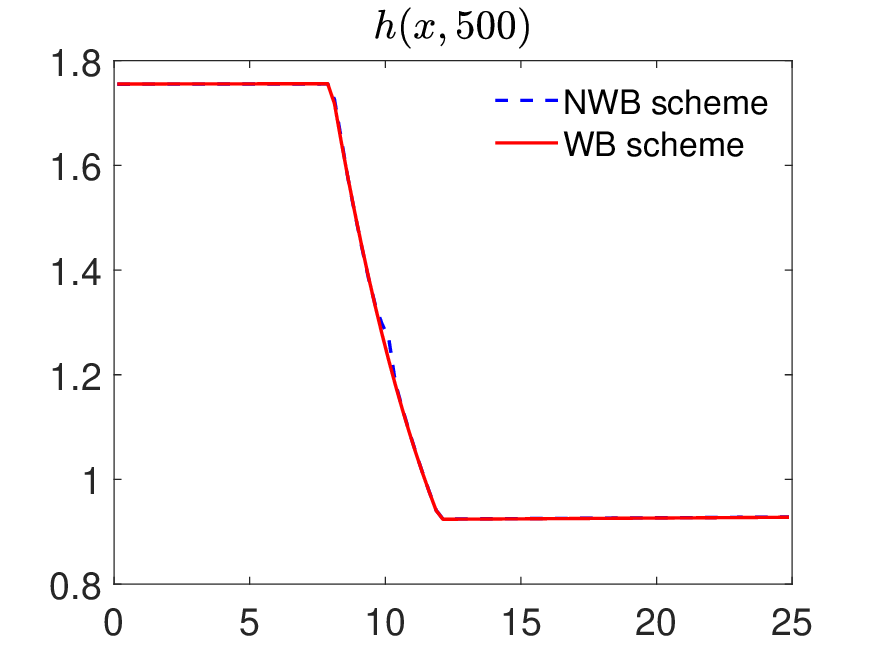}\hspace*{0.5cm}
            \includegraphics[trim=0.0cm 0.3cm 1.cm 0.1cm,clip,width=5cm]{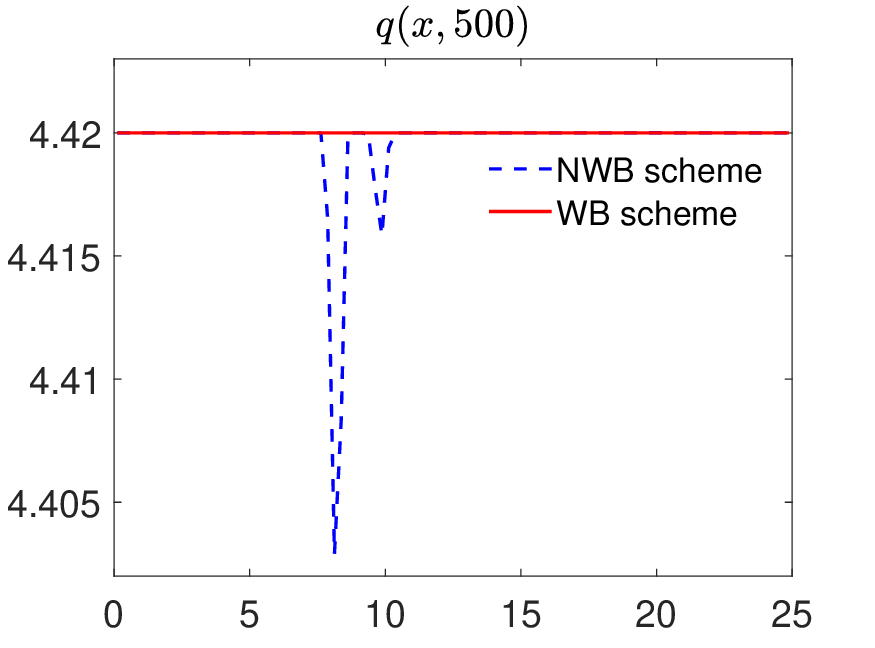}}
\vskip5pt
\centerline{\includegraphics[trim=0.0cm 0.3cm 1.cm 0.1cm,clip,width=5cm]{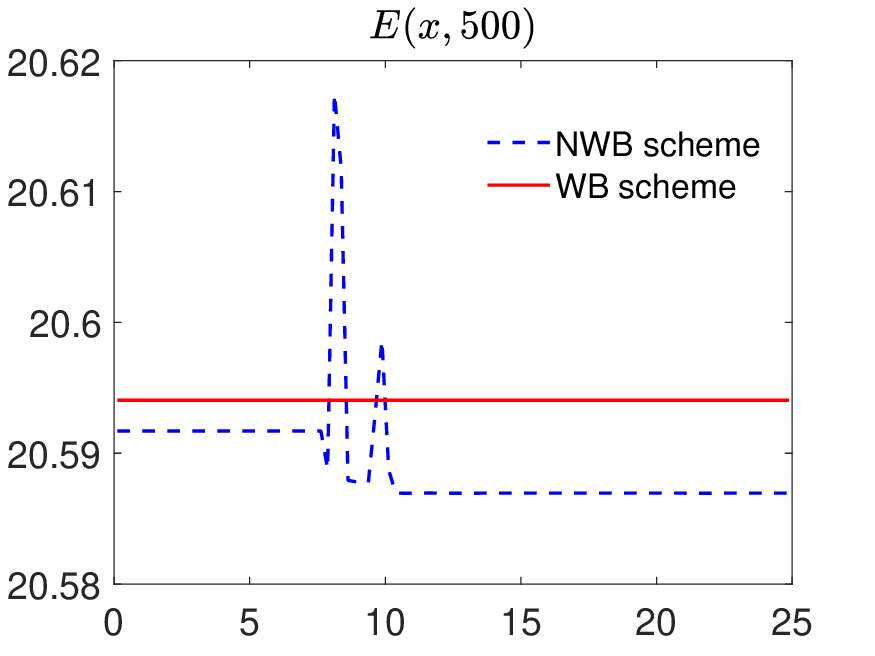}\hspace*{0.5cm}
            \includegraphics[trim=0.0cm 0.3cm 1.cm 0.1cm,clip,width=5cm]{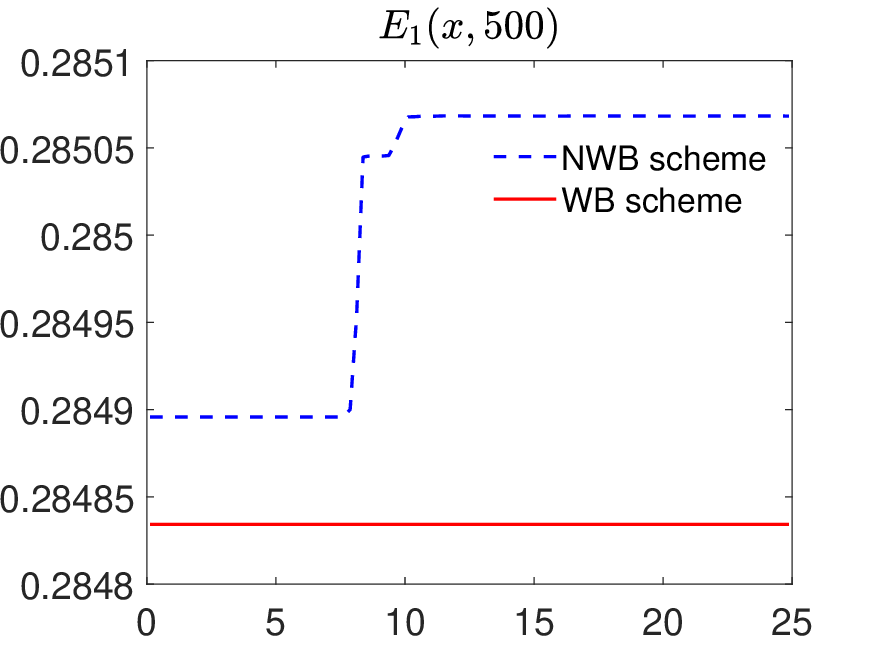}\hspace*{0.5cm}
            \includegraphics[trim=0.0cm 0.3cm 1.cm 0.1cm,clip,width=5cm]{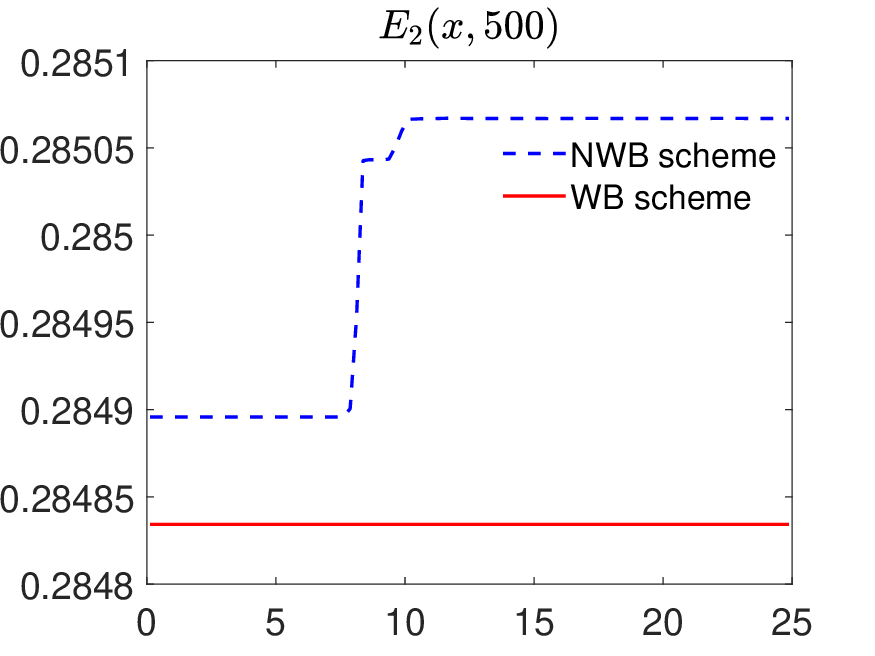}}
\vskip8pt
\centerline{\includegraphics[trim=0.0cm 0.3cm 1.cm 0.1cm,clip,width=5cm]{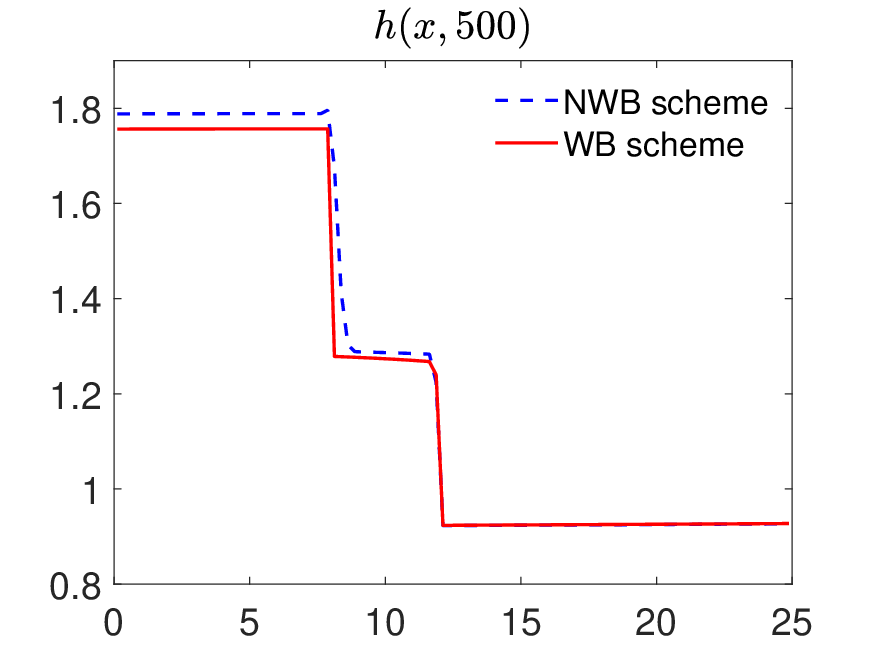}\hspace*{0.5cm}
            \includegraphics[trim=0.0cm 0.3cm 1.cm 0.1cm,clip,width=5cm]{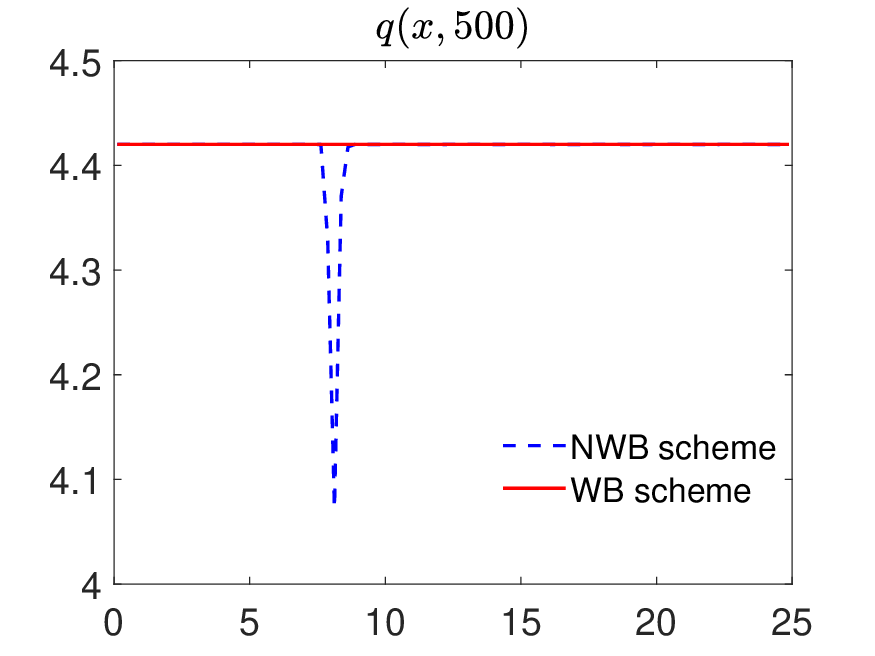}}
\vskip5pt
\centerline{\includegraphics[trim=0.0cm 0.3cm 1.cm 0.1cm,clip,width=5cm]{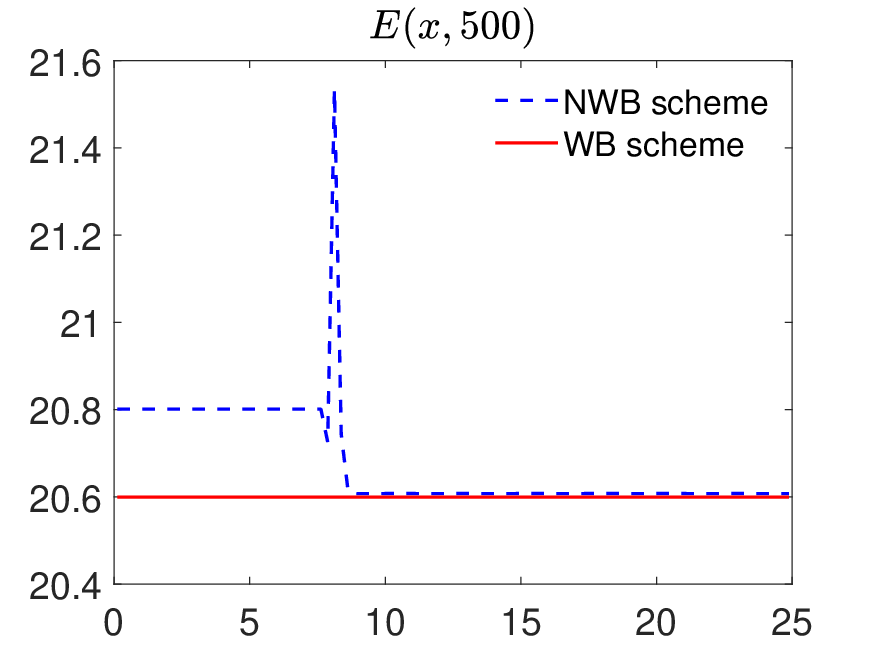}\hspace*{0.5cm}
            \includegraphics[trim=0.0cm 0.3cm 1.cm 0.1cm,clip,width=5cm]{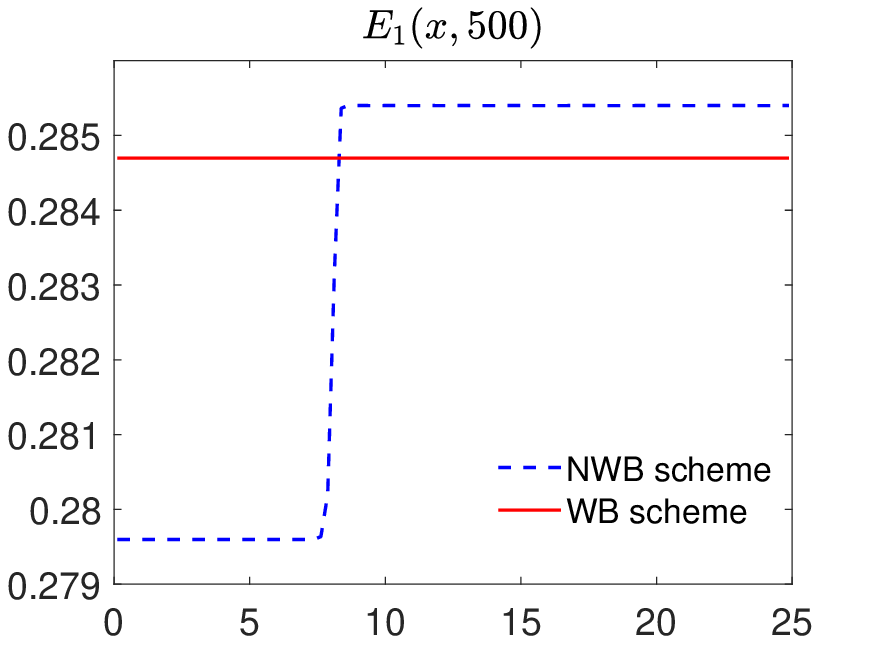}\hspace*{0.5cm}
            \includegraphics[trim=0.0cm 0.3cm 1.cm 0.1cm,clip,width=5cm]{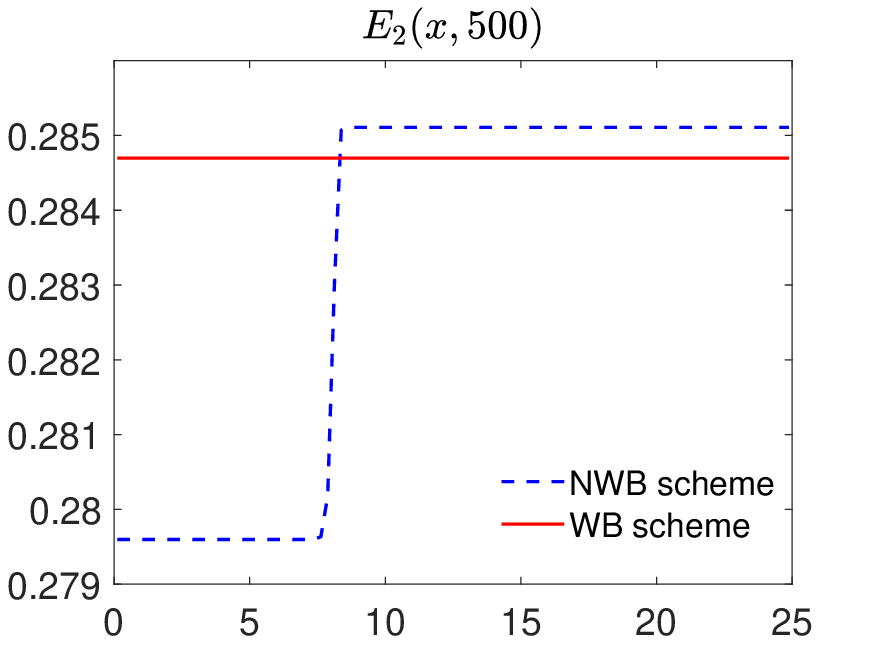}}
\caption{\sf Same as in Figures \ref{fig41} and \ref{fig42}, but for Case (c).\label{fig43}}
\end{figure}

We compute the solutions by the flux globalization based WB and NWB PCCU schemes on $100$ uniform cells until the final time $t=500$ in the
computational domain $[0,25]$ for all of the aforementioned three cases, and take the friction term $\nu=0.005$ in Cases (a) and (b), and
$\nu=0.001$ in Case (c). The obtained solutions $(h,q,E,E_1,E_2)$ for the three cases are plotted in Figures \ref{fig41}--\ref{fig43},
respectively, where it can be observed that the computed solutions converge, as expected, to the corresponding discrete steady states for
both continuous and discontinuous bottom topographies using the two schemes in Case (a); see Figure \ref{fig41}. One can also observe that
when the bottom topography is continuous, the solutions computed by the WB PCCU scheme are more accurate in Cases (b) and (c): $q$, $E$,
$E_1$, and $E_2$ computed by the WB PCCU scheme are almost constant, while the NWB PCCU solutions generate noticeable nonphysical
oscillations; see the top two rows in Figures \ref{fig42} and \ref{fig43}. When the bottom topography is discontinuous, the NWB solution
seems to be inaccurate in both Cases (b) and (c); see the bottom two rows in Figures \ref{fig42} and \ref{fig43}.

\subsubsection*{Example 2---Small Perturbations of the Steady States}
\begin{figure}[ht!]
\centerline{\includegraphics[trim=0.6cm 0.3cm 1.1cm 0.1cm,clip,width=5cm]{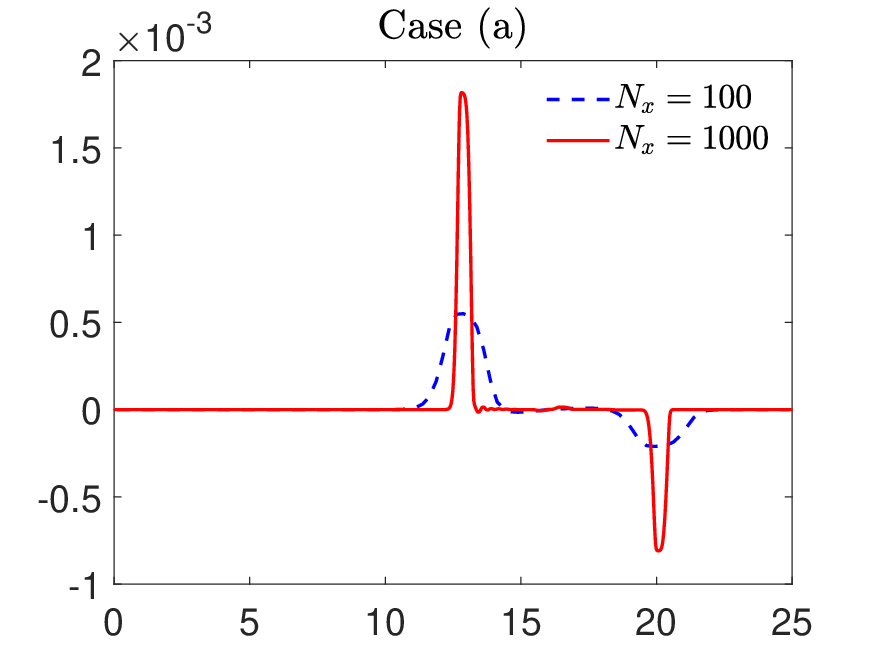}\hspace*{0.5cm}
            \includegraphics[trim=0.6cm 0.3cm 1.1cm 0.1cm,clip,width=5cm]{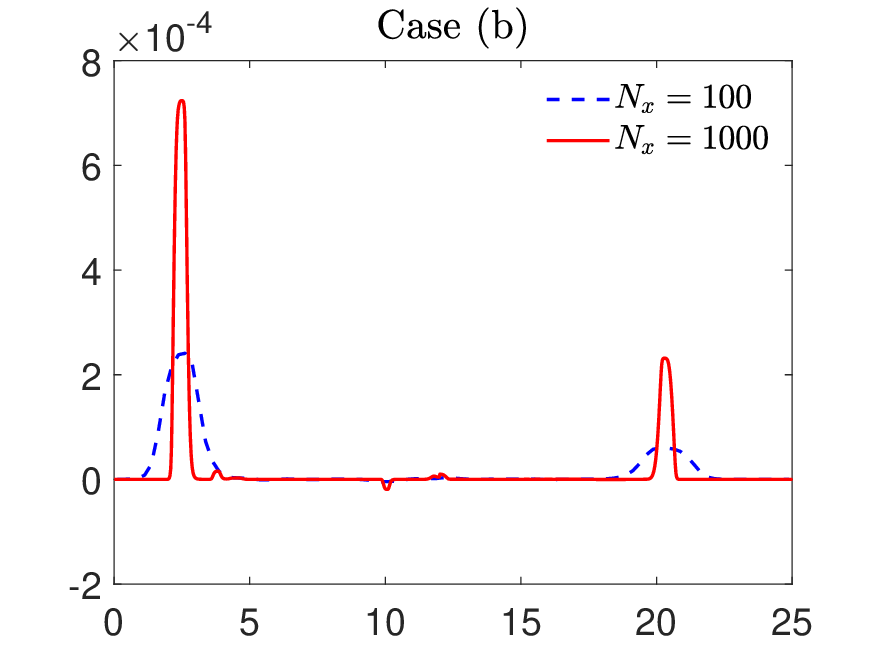}\hspace*{0.5cm}
            \includegraphics[trim=0.6cm 0.3cm 1.1cm 0.1cm,clip,width=5cm]{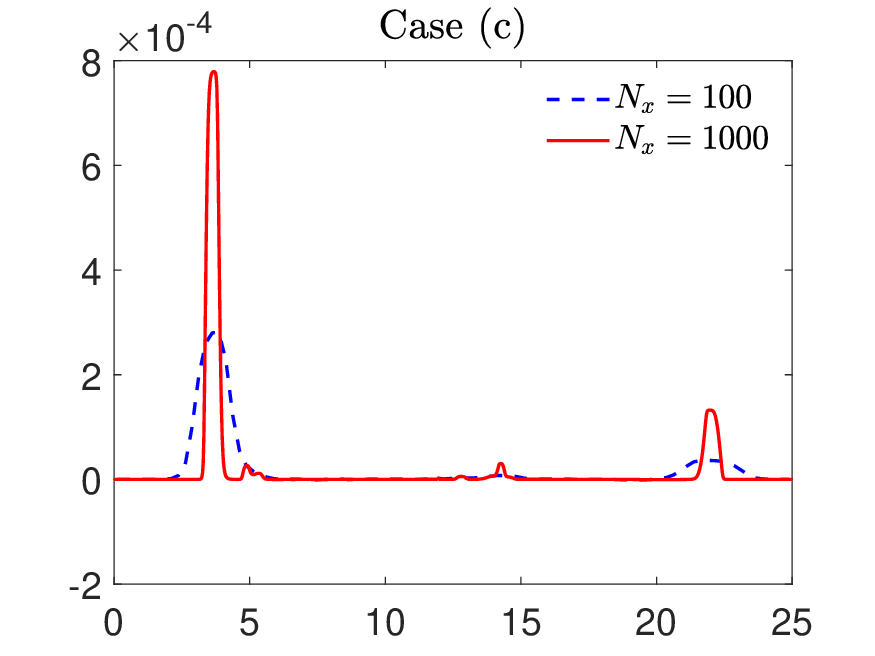}}
\caption{\sf Example 2 (small perturbation): The difference $h(x,t)-h_{\rm eq}(x)$ computed using $100$ and $1000$ uniform cells at times
$t=0.8$ in Case (a) and $t=2$ in Cases (b) and (c).\label{fig44}}
\end{figure}

In the second example, we test the ability of the WB PCCU scheme to capture small perturbations of the discrete equilibria computed in
Example 1, in the case of discontinuous bottom topography given by \eref{4.1}. The initial data are obtained by adding a small positive
number $0.001$ in the interval $x\in[5.75,6.25]$ to the steady-state water depth $h_{\rm eq}(x)$ obtained at the end of the WB PCCU
simulations reported in Example 1. We compute the solutions until the final times $t=0.8$ in Case (a) and $t=2$ in Cases (b) and (c) using
both $100$ and $1000$ uniform cells in the computational domain $[0,25]$. In Figure \ref{fig44}, where we plot the difference
$h-h_{\rm eq}$, one can observe that no spurious oscillations are produced and the propagating perturbations are well captured by the WB
PCCU scheme.

\subsubsection*{Example 3---Dam-Break Problem}
In this example, we  consider the discontinuous initial data corresponding to an extreme velocity profile (taken from
\cite{KR_2020,KP_SWM}), which correspond to a dam break:
\begin{equation*}
\begin{aligned}
&h(x,0)=\left\{\begin{aligned}
&5,&&\mbox{if}~x<0,\\
&1,&&\mbox{if}~x>0,
\end{aligned}\right.
\qquad u(x,0)\equiv0.25,\\
&\alpha_1(x,0)\equiv-0.25,\quad\alpha_N(x,0)\equiv0.25,\quad\alpha_i(x,0)\equiv0,\quad i=2,\ldots,N-1.
\end{aligned}
\end{equation*}
We take $N=8$, $g=1$, and a flat bottom topography $Z(x)\equiv0$.

\begin{figure}[ht!]
\centerline{\includegraphics[trim=0.6cm 0.3cm 1.0cm 0.1cm,clip,width=5cm]{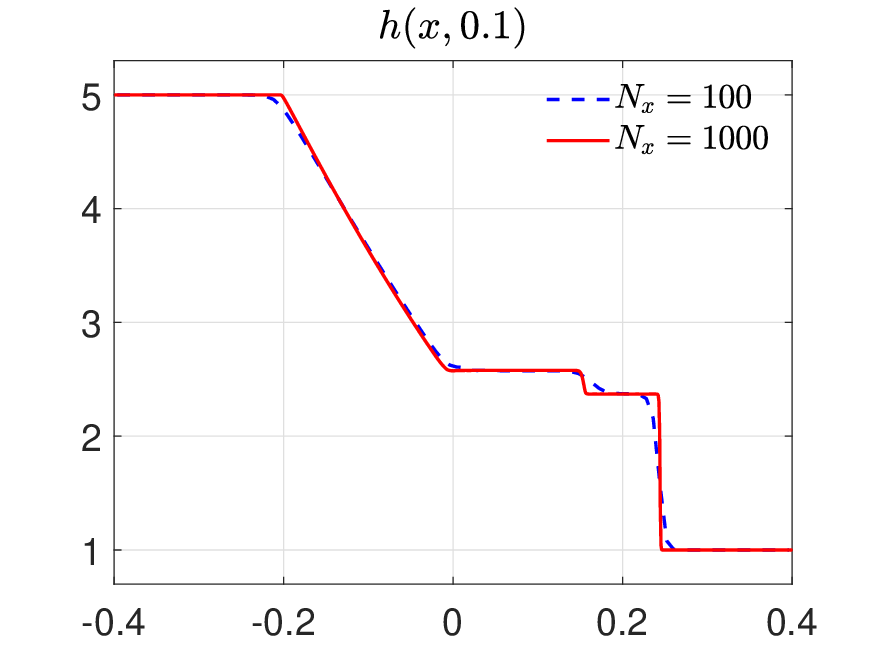}\hspace*{0.5cm}
            \includegraphics[trim=0.6cm 0.3cm 1.0cm 0.1cm,clip,width=5cm]{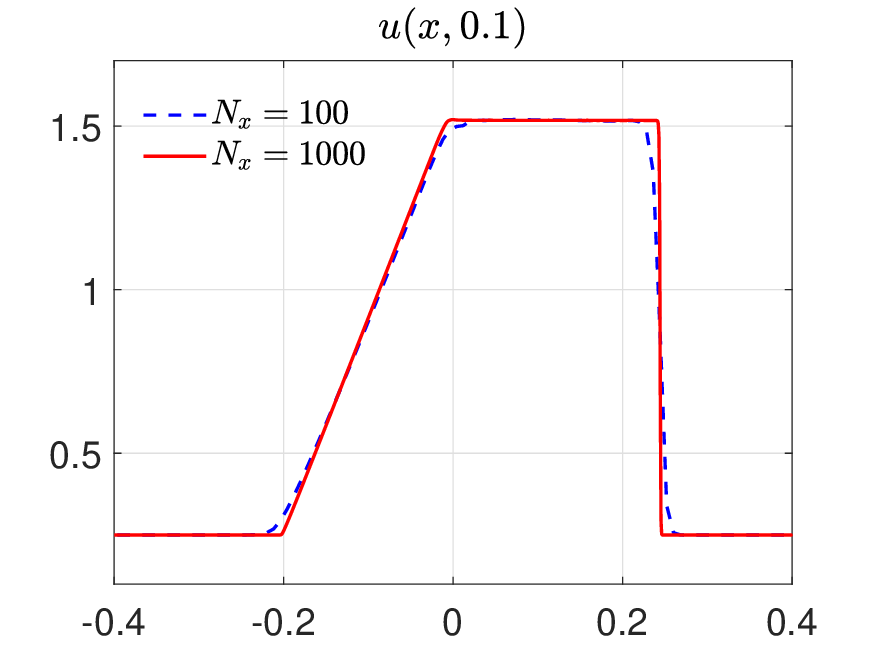}}
\vskip5pt
\centerline{\includegraphics[trim=0.6cm 0.3cm 1.0cm 0.1cm,clip,width=5cm]{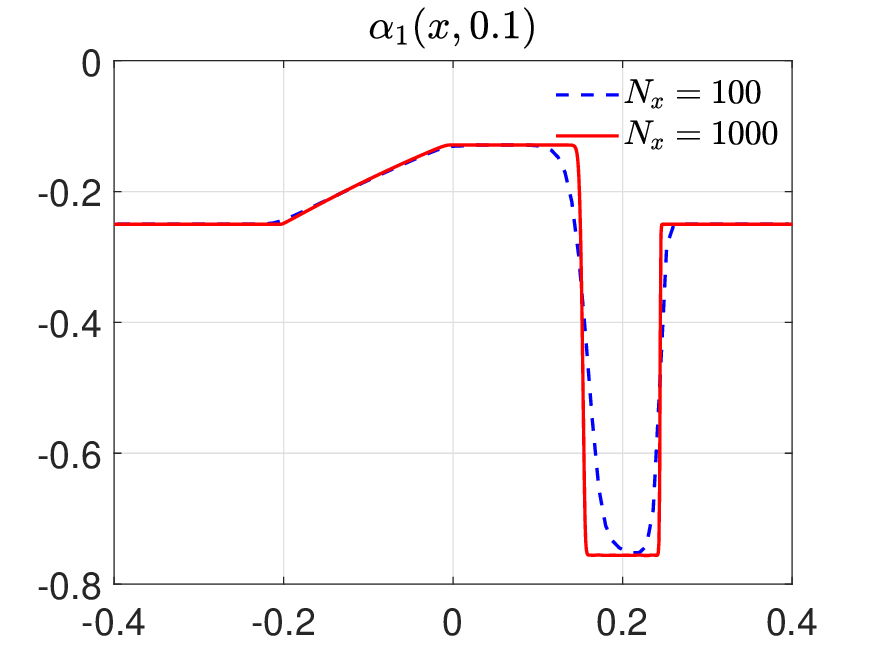}\hspace*{0.5cm}
            \includegraphics[trim=0.6cm 0.3cm 1.0cm 0.1cm,clip,width=5cm]{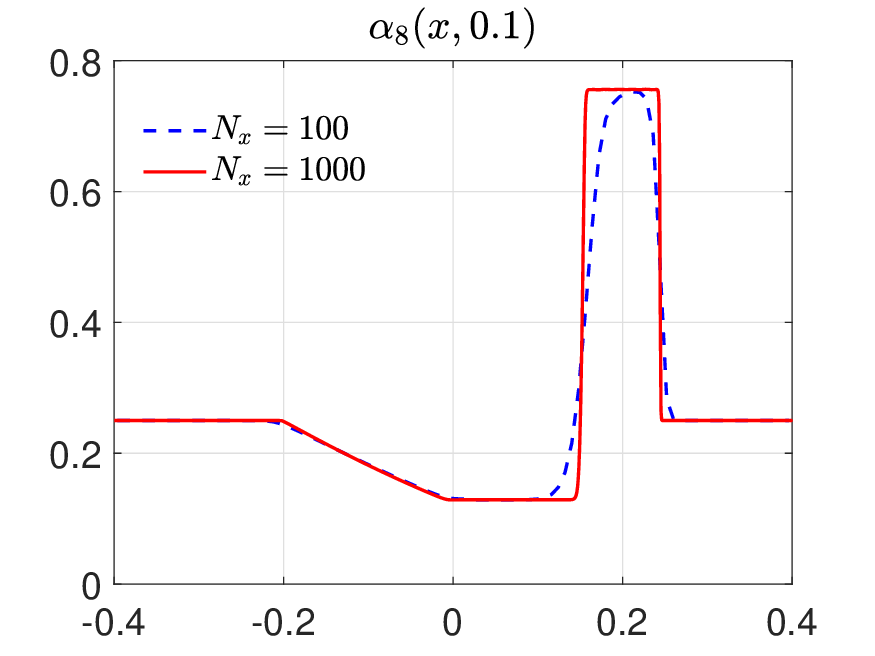}}
\caption{\sf Example 3: $h$, $u$, $\alpha_1$, and $\alpha_8$ computed by the WB PCCU scheme. $\nu=0$.\label{fig45}}
\end{figure}

\begin{figure}[ht!]
\centerline{\includegraphics[trim=0.6cm 0.3cm 1.0cm 0.1cm,clip,width=5cm]{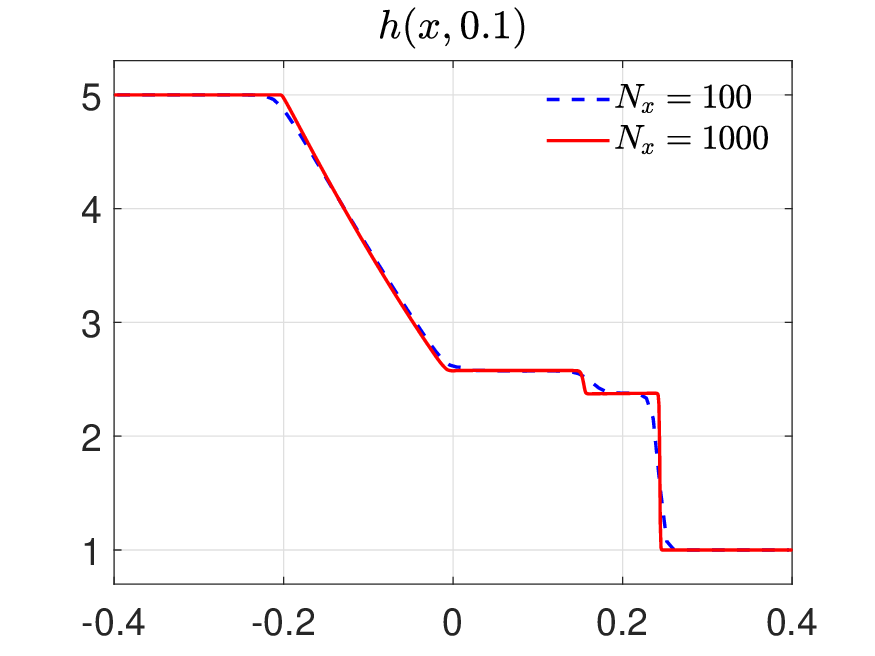}\hspace*{0.5cm}
            \includegraphics[trim=0.6cm 0.3cm 1.0cm 0.1cm,clip,width=5cm]{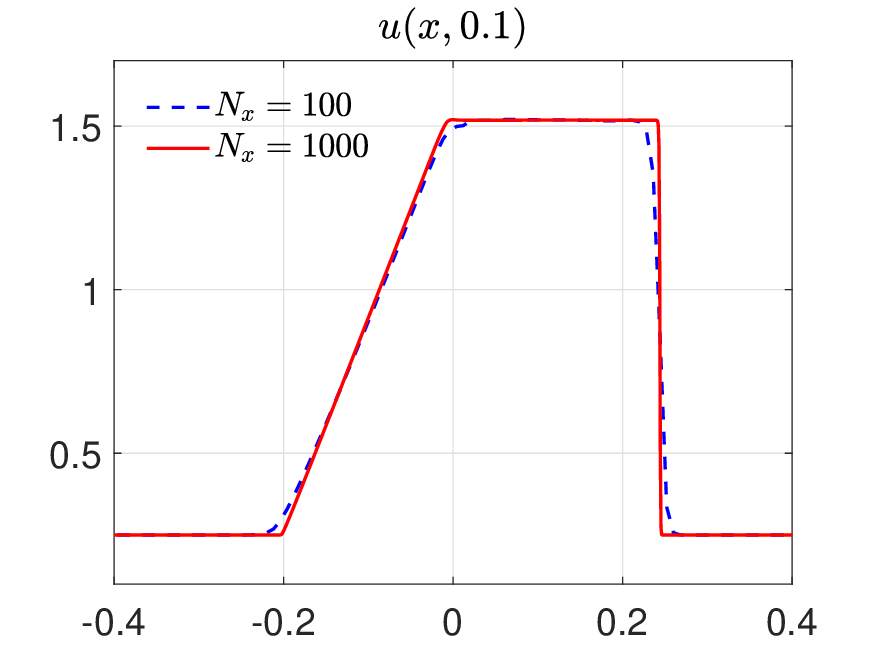}}
\vskip5pt
\centerline{\includegraphics[trim=0.6cm 0.3cm 1.0cm 0.1cm,clip,width=5cm]{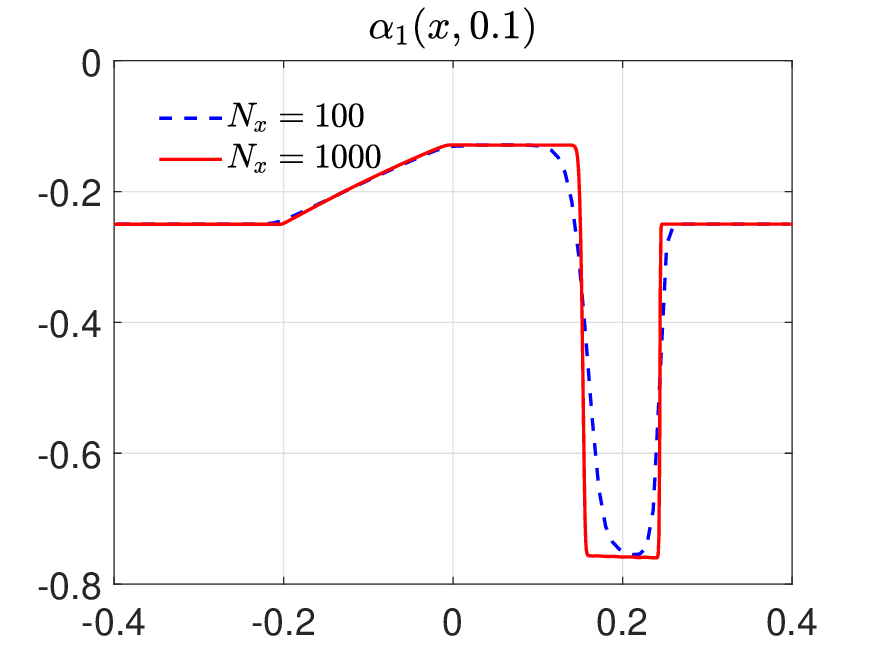}\hspace*{0.5cm}
            \includegraphics[trim=0.6cm 0.3cm 1.0cm 0.1cm,clip,width=5cm]{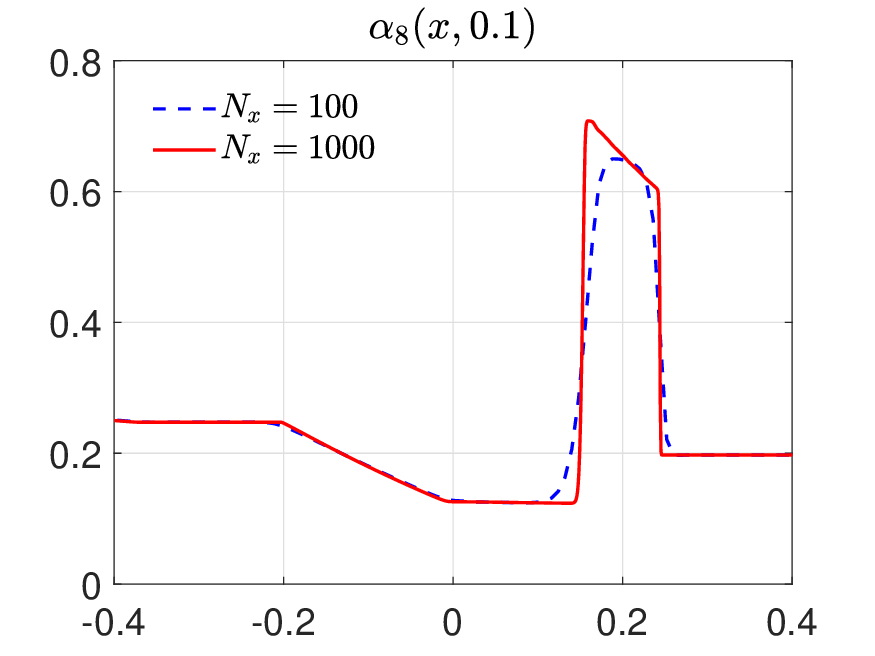}}
\caption{\sf Example 3: Same as in Figure \ref{fig45}, but with the friction coefficient $\nu=0.001$.\label{fig46}}
\end{figure}

We compute the numerical solutions using the WB PCCU scheme until the finial time $t=0.1$ using both $100$ and $1000$ uniform cells in the
spatial domain $[-0.4,0.4]$. The obtained $h$, $u$, $\alpha_1$, and $\alpha_8$ are plotted in Figures \ref{fig45} and \ref{fig46} for both
$\nu=0$ (frictionless case) and $\nu=0.001$. One can clearly see that the scheme captures the shock and rarefaction waves very well. It
should be also noted that the WB PCCU scheme outperforms the scheme from \cite{KP_SWM}, which generates small nonphysical waves near $x=0$.
We also stress that in \cite{KP_SWM}, only the frictionless case was considered. The friction case includes a notably different solution for
the last moment $\alpha_8$ in Figure \ref{fig46} (bottom right).

\subsubsection*{Example 4--- Dam-Break Problem With Square Root Velocity Profile}
In the final example also taken from \cite{KR_2020,KP_SWM}, we consider the same initial set-up as in Example 3, but with a more complex
velocity profile resulting from the expansion of a square-root shape with $N=8$ coefficients, such that $u(x,0)\equiv1$ and
\begin{equation*}
\begin{aligned}
&\alpha_1(x,0)\equiv-\frac{3}{5},&&\alpha_2(x,0)\equiv-\frac{1}{7},&&\alpha_3(x,0)\equiv-\frac{1}{15},&&\alpha_4(x,0)\equiv-\frac{3}{77},\\
&\alpha_5(x,0)\equiv-\frac{1}{39},&&\alpha_6(x,0)\equiv-\frac{1}{55},&&\alpha_7(x,0)\equiv-\frac{3}{221},&&\alpha_8(x,0)\equiv-\frac{1}{95}.
\end{aligned}
\end{equation*}

\begin{figure}[ht!]
\centerline{\includegraphics[trim=0.6cm 0.3cm 1.0cm 0.1cm,clip,width=5cm]{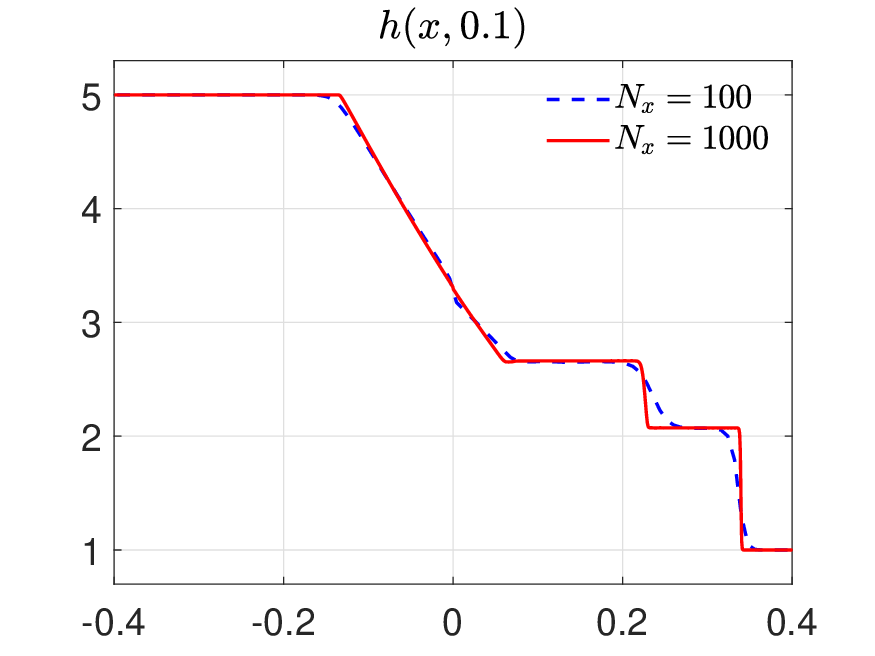}\hspace*{0.5cm}
            \includegraphics[trim=0.6cm 0.3cm 1.0cm 0.1cm,clip,width=5cm]{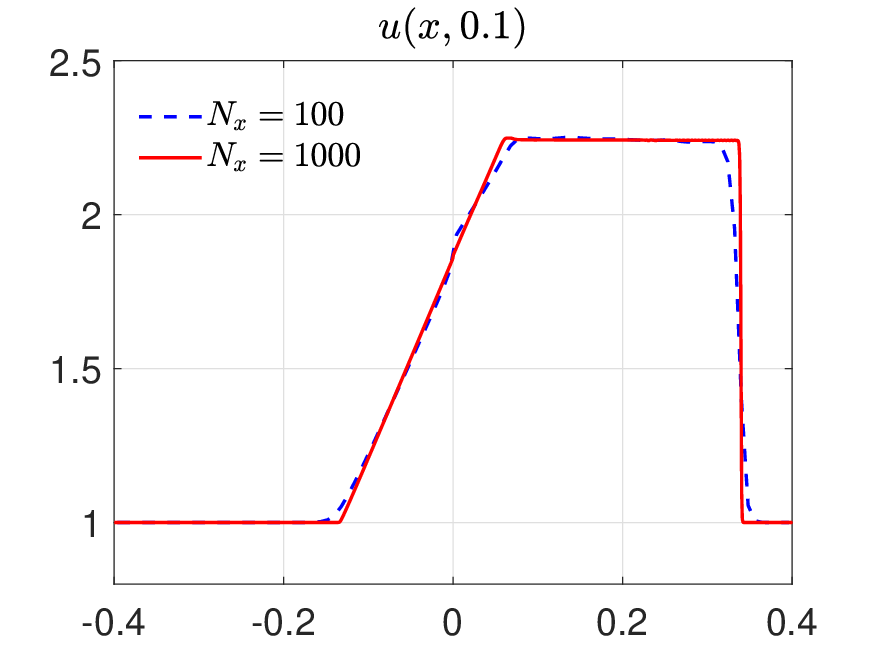}}
\vskip5pt
\centerline{\includegraphics[trim=0.6cm 0.3cm 1.0cm 0.1cm,clip,width=5cm]{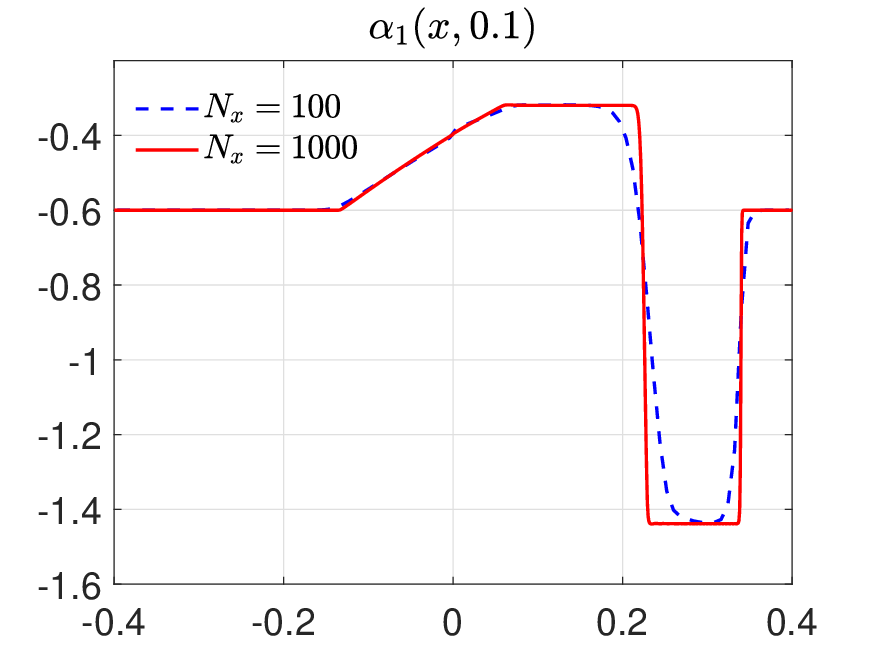}\hspace*{0.5cm}
            \includegraphics[trim=0.6cm 0.3cm 1.0cm 0.1cm,clip,width=5cm]{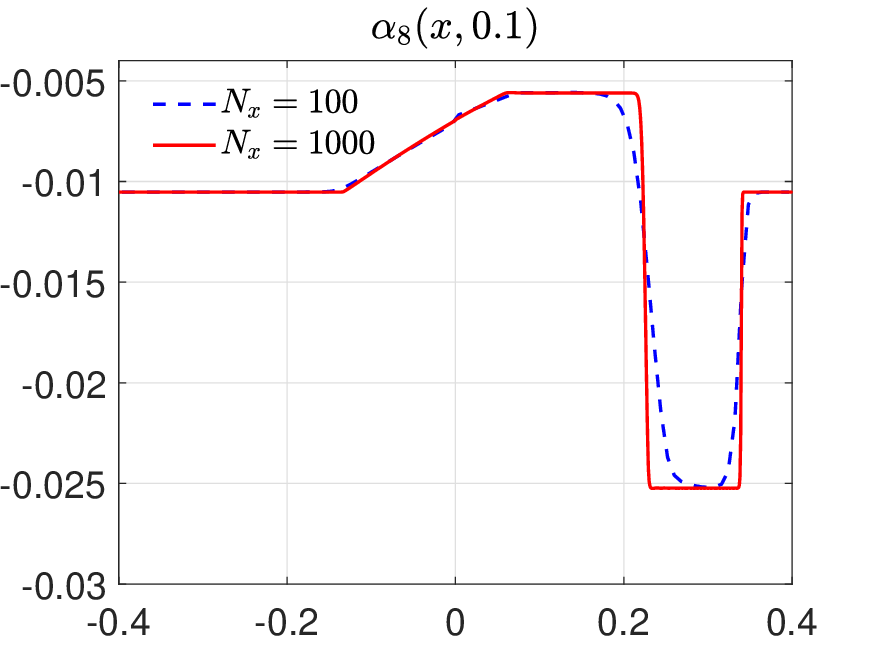}}
\caption{\sf  Example 4: $h$, $u$, $\alpha_1$, and $\alpha_8$ computed by the WB PCCU scheme. $\nu=0$.\label{fig47}}
\end{figure}

\begin{figure}[ht!]
\centerline{\includegraphics[trim=0.6cm 0.3cm 1.0cm 0.1cm,clip,width=5cm]{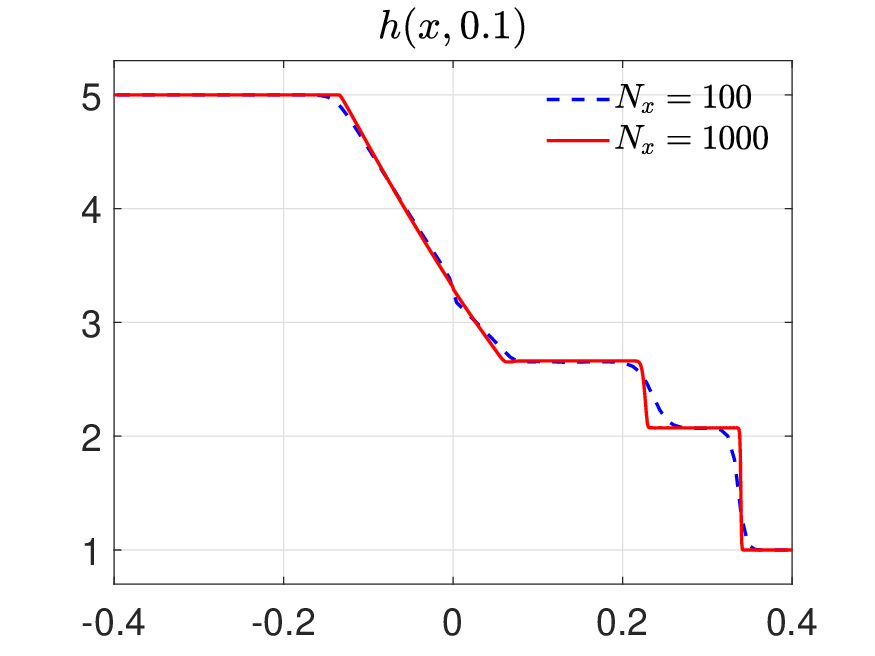}\hspace*{0.5cm}
            \includegraphics[trim=0.6cm 0.3cm 1.0cm 0.1cm,clip,width=5cm]{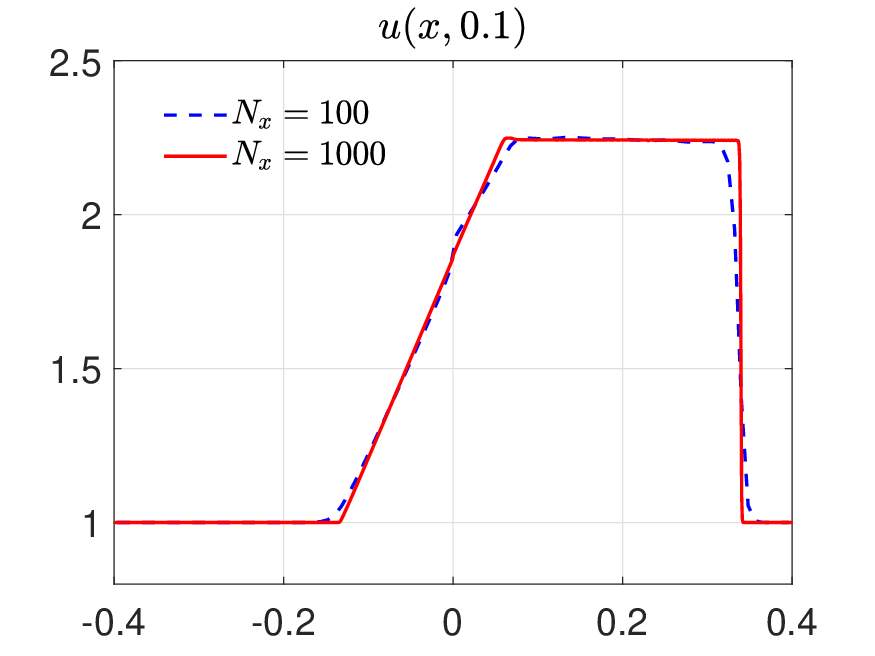}}
\vskip5pt
\centerline{\includegraphics[trim=0.6cm 0.3cm 1.0cm 0.1cm,clip,width=5cm]{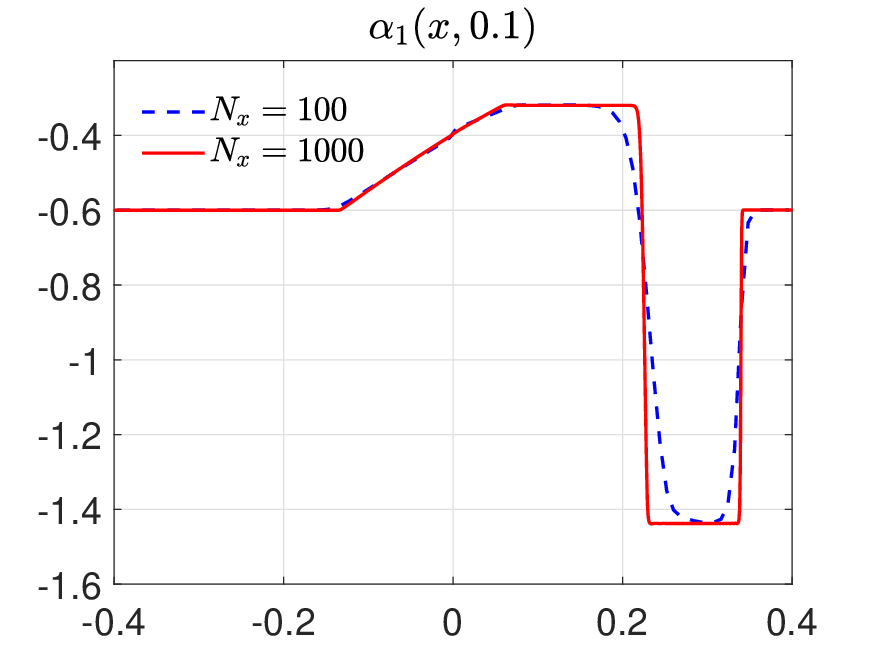}\hspace*{0.5cm}
            \includegraphics[trim=0.6cm 0.3cm 1.0cm 0.1cm,clip,width=5cm]{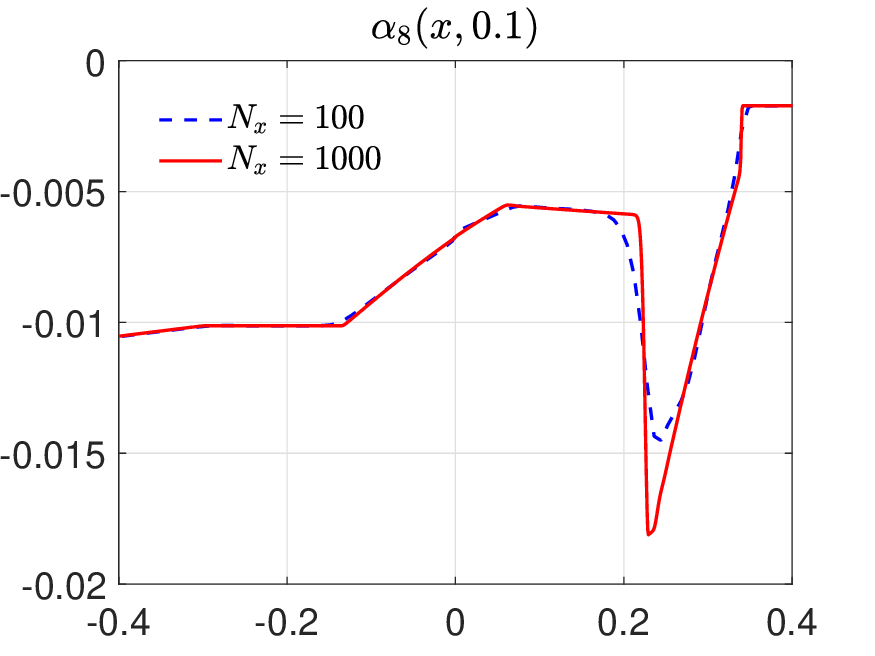}}
\caption{\sf  Example 4: Same as in Figure \ref{fig47}, but with the friction coefficient $\nu=0.001$.\label{fig48}}
\end{figure}

We compute the numerical solutions using the WB PCCU scheme until the finial time $t=0.1$ using both $100$ and $1000$ uniform cells in the
spatial domain $[-0.4,0.4]$. The obtained $h$, $u$, $\alpha_1$, and $\alpha_8$ are plotted in Figures \ref{fig47} and \ref{fig48} for both
$\nu=0$ (frictionless case) and $\nu=0.001$. As in Example 3, one can clearly see that the proposed WB PCCU scheme accurately captures the
solution. One can also observe that the computed solution does not contain any nonphysical waves near $x=0$, which were generated by the
scheme from \cite{KP_SWM}.

%
%

\section{Conclusion and future work}
We have developed flux globalization based second-order well-balanced path-conservative central-upwind schemes for the hyperbolic shallow
water linearized moment equations. The proposed schemes combine flux globalization and path-conservative formulations tailored to
nonconservative products for a system of equations. The conducted numerical experiments confirm that the methods are accurate and
well-balanced. The obtained results also demonstrate the benefit of including higher-order moments to resolve vertical variations. Future
work includes extensions to more complex shallow water moment models, two-dimensional models, and the inclusion of more physical effects
like bedload transport, sediments, or different friction models.

\subsubsection*{Acknowledgment}
The work of Y. Cao was supported in part by the Guangdong Basic and Applied Basic Research Foundation (No. 2025A1515012249). The work of
Q. Huang was supported by the Deutsche Forschungsgemeinschaft (DFG, German Research Foundation) – SPP 2410 Hyperbolic Balance Laws in Fluid
Mechanics: Complexity, Scales, Randomness (CoScaRa). The work of A. Kurganov was supported in part by the NSFC grant 12171226 and W2431004.
The work of Y. Liu was supported in part by the SNFS grants 200020\_204917 and UZH Posdoc Grant, 2024/Verf\"{u}gung Nr.FK-24-110. The
publication is part of the project \textit{HiWAVE} with file number VI.Vidi.233.066 of the research \textit{ENW Vidi} programme, funded by
the \textit{Dutch Research Council (NWO)} under the grant \url{https://doi.org/10.61686/CBVAB59929}.

\appendix

\section{Generalized Minmod Piecewise Linear Reconstruction}\label{appxA}
In this appendix, we describe a second-order piecewise linear reconstruction based on the generalized minmod limiter; see, e.g.,
\cite{LN,NT,Swe}.

We assume that the values $\varphi_k$ (either the cell averages or point values) of a certain function $\varphi$ at the cell centers $x=x_k$
are available. Then a second-order non-oscillatory piecewise linear reconstruction of $\varphi$ reads as
\begin{equation*}
\widetilde{\varphi}(x)=\varphi_k+(\varphi_x)_k(x-x_k),\quad x\in C_k,
\end{equation*}
where the slopes $(\varphi_x)_k$ have to be evaluated using a nonlinear limiter. In our numerical experiments, we have used the generalized
minmod limiter (see \cite{LN,NT,Swe}):
\begin{equation*}
(\varphi_x)_k={\rm minmod}\left(\theta\,\frac{\varphi_{k+1}-\varphi_k}{\dx},\,\frac{\varphi_{k+1}-\varphi_{k-1}}{2\dx},\,
\theta\,\frac{\varphi_k-\varphi_{k-1}}{\dx}\right),\quad\theta\in[1,2],
\end{equation*}
where the minmod function is defined by
\begin{equation*}
\mbox{minmod}(c_1,c_2,\ldots)=\left\{\begin{aligned}
&\min(c_1,c_2,\ldots)&&\mbox{if}~c_i>0,~\forall i,\\
&\max(c_1,c_2,\ldots)&&\mbox{if}~c_i<0,~\forall i,\\
&0&&\mbox{otherwise},
\end{aligned}\right.
\end{equation*}
and the parameter $\theta$ is used to control the amount of numerical dissipation present in the resulting scheme: larger values of $\theta$
correspond to sharper but, in general, more oscillatory reconstructions.

\bibliography{reference}
\bibliographystyle{siam}
\end{document}